\input amstex

\loadeufm
\loadmsbm
\loadeufm

\documentstyle{amsppt}
\input amstex
\catcode `\@=11
\def\logo@{}
\catcode `\@=11
\magnification \magstep1
\NoRunningHeads
\NoBlackBoxes
\TagsOnLeft

\def \={\ = \ }
\def \+{\ +\ }
\def \-{\ - \ }

\def \b|{\big |}

\def \g1{\Gamma_1}

\def \nfp{\demo\nofrills{Proof:\usualspace\usualspace }}

\def\rarr#1#2{\smash{\mathop{\hbox to .5in{\rightarrowfill}}
         \limits^{\scriptstyle#1}_{\scriptstyle#2}}}

\def\larr#1#2{\smash{\mathop{\hbox to .5in{\leftarrowfill}}
          \limits^{\scriptstyle#1}_{\scriptstyle#2}}}

\def\swarr#1#2 {\llap{$\scriptstyle #1$}  \swarrow
        \vcenter to .5in{}\rlap{$\scriptstyle #2$}}

\topmatter
\title
\centerline{ A canonical  Makanin-Razborov diagram and a pseudo  }
\centerline{ topology for sets of tuples in free groups,  }
\centerline {semigroups, associative
algebras and Lie algebras I}
\endtitle
\author
\centerline{
Z. Sela${}^{1,2}$}
\endauthor

\dedicatory{To the memory of Eliyahu Rips}
\enddedicatory

\footnote""{${}^1$Hebrew University, Jerusalem 91904, Israel.}
\footnote""{${}^2$Partially supported by an Israel academy of sciences fellowship.}
\abstract\nofrills{}
The JSJ decomposition and the Makanin-Razborov diagram were proved to be  essential in studying varieties over free
groups, semigroups and associative algebras. In this paper we suggest a unified conceptual approach to
the applicability of these structures over all these algebraic categories.

With a variety over each of these algebraic categories we naturally associate a set of tuples in a free group.
Then we show how to associate a Makanin-Razborov diagram with any
set of tuples over a free group.

Furthermore, in case the MR diagram that is associated with a set of tuples is single ended, we prove that there is a canonical Makanin-Razborov diagram that
can be associated with such a set. This canonical diagram is a main key in studying varieties over free semigroups,
associative algebras and Lie algebras, and encodes the global structure of these varieties. It enables us to define a (pseudo) closure of a set of tuples
over each of the algebraic objects, associate a rank with it (analogous to Shelah and Lascar ranks), and  over free groups the closure provides a canonical envelope
that is essential in  studying  the structure and the properties of definable sets.
\endabstract
\endtopmatter

\document

\baselineskip 12pt
Varieties over free groups were studied by Makanin [Ma] and Razborov [Ra1],[Ra2]. In [Se1] we presented a geometric approach to study these varieties
that uses the canonical JSJ decomposition of groups to construct a Makanin-Razborov diagram of a variety, that encodes the set of solutions to a
system of equations over a free group. The whole construction was later generalized to varieties over hyperbolic groups
([Se6],[Re-We]) and over other groups. The construction of the MR diagram for varieties is canonical if we fix a generating set of the involved limit groups,
but is not canonical
when these generating sets are not fixed  as was shown by G. Berk [Be].

Varieties over free semigroups were first studied by A. Jez [Je] using recompression methods. In [Se9] we constructed a Makanin-Razborov diagram
that encodes the solutions to a system of equations over a free semigroup. The construction uses a compactness argument, so the constructed
diagram is not canonical.

In an ongoing work with A. Atkarskaya we study the structure of varieties over a free associative algebra [At-Se]. The key to our analysis is  the
Makanin-Razborov diagram that was constructed over free semigroups, that encodes the top homogeneous parts of the set of solutions to a system
of equations over the associative algebra.

Free Lie algebras are naturally embedded into free associative algebras, so every variety over a free Lie algebra can be viewed as the intersection of
a variety over a free associative algebra with the appropriate embedding of a product of the free Lie algebra into a product of the associative algebra.

The study of all these varieties uses a Makanin-Razborov diagram in an essential way, and through it low dimensional topology is central in the
study of varieties over free objects in the corresponding categories. In this paper we present a unification of  the constructions of the Makanin-Razborov
diagrams that were constructed for varieties over free objects in the various categories.

Instead of associating a Makanin-Razborov diagram with a variety we construct an MR diagram for every set of $\ell$-tuples in a free group $F_k$. The diagram
that is associated with a set is in general not canonical, as it is based on a compactness argument, but for sets for which there exists a single ended diagram
that can be associated with them, we prove that it is possible to replace the diagram with a canonical one.

A $single$ $ended$ MR diagram is a diagram for which all the abelian decompositions along its levels contain no free products. First order
formulas that have a similar property were termed $minimal$ $rank$ in [Se3],[Se4]. In the work on
Tarski  problem it is shown that the analysis of minimal rank formulas is considerably easier than the analysis of general formulas.

The first order theory of a non-abelian free group is stable, but it is not superstable [Po1]. Indeed,
the minimal rank formulas are precisely the superstable part of the theory of a free group. Superstable definable sets satisfy a certain descending chain
condition, that enabled Shelah, Lascar and others to associate a countable ordinal with such a definable set.

We associate a (geometric) complexity with a single ended diagram. The complexities of single ended diagrams are well-ordered. Hence, with a set of tuples
for which it is possible to associate a single ended MR diagram, it is possible to associate a minimal complexity MR diagram. In theorem 1.16 we prove that
this minimal complexity diagram is essentially unique (up to an appropriate isomorphism), and hence an invariant of the set of tuples.

In the second section we use this canonical diagram  to define a pseudo closure of a set of tuples $S$ in a free group, as the set of tuples that extend
to homomorphisms that factor through
the canonical diagram of $S$. Every (minimal rank) variety is pseudo closed, so the notion of pseudo closed considerably refine the Zariski topology over a free group.
Like the Zariski topology, pseudo closed sets are a Noetherian family. i.e., every sequence of pseudo closed sets
with strictly decreasing complexities terminates after a finite step.

The canonical MR diagram and the pseudo closre enable us to associate a $rank$ with a set of tuples. The rank is a countable ordinal and is defined
in a similar way to the Shelah or Lascar ranks in superstable theories for definable sets and types. Note that our definition of the rank of a set involves no
model theory.

The pseudo closure of a set enables us to associate a canonical $envelope$ with
a minimal rank definable set (theorem 2.5). (Non-canonical) envelopes of definable sets were introduced in [Se8] and were previously used to classify imaginaries,
and (explicitly) to prove the stability of the theory in [Se7].

In the 3rd section we show how to naturally associate a set of tuples in a free group with a set of tuples in a free semigroup. With the results of [Se9],
 this enables one to transfer all the constructions over free groups to free semigroups.

In the last section we explain how to associate a set of tuples in a free semigroup with a set of tuples in a free associative algebra. The set of tuples in a semigroup
is really associated with the top homogeneous parts of the tuples in the associative algebra. In an ongoing work with A. Atkarskaya [At-Se] we show how to use the canonical
diagram that we associate with the set of tuples in a variety over a free associative algebra to study the structure of the variety. Finally, since a free Lie algebra embeds
in a free associative algebra, we can associate a canonical diagram with a set of tuples over a free Lie algebra, and the canonical diagram enables one to study
the structure of varieties over Lie algebras.

\vglue 1.5pc
\centerline{\bf{\S1.  A canonical Makanin-Razborov diagram for a subset}}
\centerline{\bf{in a power of a free group}}
\medskip

Let $F_k$ be a non-abelian free group, and let $G$ be a f.g.\ group. In section 5 in [Se1] we have associated a canonical Makanin-Razborov diagram
with $Hom(G,F_k)$, where $Hom(G,F_k)$ is equivalent to a variety over a free group.
With $Hom(G,F_k)$ we first associated  the canonical finite collection of maximal limit quotients of $G$, $L_1,\ldots,L_s$.
With each limit quotient, $L_i$,  we associated the homomorphisms that factor through it, i.e., $h:G \to F_k$ that can be written
as $h=f \circ \eta_i$, where $\eta_i$ is the canonical quotient map: $G \to L_i$, and $f: L_i \to F_k$. The collection of homomorphisms
from $Hom(G,F_k)$ that factor through $L_i$ is the irreducible subvariety of $Hom(G,F_k)$ with a dual group $L_i$, i.e., the irreducible
subvariety $Hom(L_i,F_k)$.

With $L_i$ we associated its JSJ decomposition, and with the JSJ decomposition we associated its modular groups. We used the modular groups
to shorten homomorphisms, i.e., we looked at shortened homomorphisms $f \circ \varphi_f$, where $f:L_i \to F_k$ and $\varphi_f$ is a shortening
automorphism from the modular group of $L_i$ (that depends on $f$). Finally, the collection of shortened homomorphisms has a strictly smaller
Zariski closure, i.e., the collection is a subset in a finite union of varieties $Hom(MSQ_1,F_k),\ldots,Hom(MSQ_t,F_k)$, where each $MSQ_j$ is
a limit quotient of the limit group $L_i$.

We continued iteratively, and by the d.c.c.\ for limit groups (the Noetherianity of free groups), the construction terminates after
finitely many steps, and the final output is the Makanin-Razborov (MR) diagram that is associated with a variety. Note that the terminal limit group
in each path (resolution) in the diagram is a free group.

In this chapter our goal is to associate an MR diagram with a general set in a power of a free group, and not just with a variety. Of course, it is possible
to associate a diagram with the set by starting with the Zariski closure of the set, and with the Zariski closure we already associated a canonical
diagram. However, the diagram that we associate with a set depends on the set and not just on its Zariski closure. In fact, it is the most efficient
MR diagram that can be associated with the set, i.e., the diagram with the lowest possible (natural) complexity. Such a diagram is crucial in studying
varieties over free semigroups and associative algebras, and has applications in model theory as well.

\smallskip

Let $\ell$ be a fixed positive integer, and let $U \subset F_k^{\ell}$ be a set of (ordered) $\ell$-tuples from $F_k$. We fix
a free basis of $F_{\ell+k}$ and a free basis of $F_k$. With  each $\ell$-tuple in $U$ we
associate a homomorphism $h:F_{k+\ell} \to F_k$, such that the first $k$ basis elements in $F_{\ell+k}$ are mapped to the fixed basis elements
of $F_k$, and the last $\ell$ basis elements in $F_{\ell+k}$ are mapped to the $\ell$-tuple in $U$.

To construct an MR diagram for the set $U$ we look at all the sequences of homomorphisms:
$\{h_n\}_{n=1}^{\infty}:F_{k+\ell} \to F_k$ where $h_n$ are distinct homomorphisms that are associated with $\ell$-tuples in $U$.

We fix such a sequence $\{h_n\}$. Since the homomorphisms are distinct, the sequence is divergent, i.e., the maximal lengths of the images of the fixed
basis of $F_{k+\ell}$ under $h_n$ approaches $\infty$ when $n$ approaches $\infty$.
We can extract a subsequence of $\{h_n\}$ for which the actions of $F_{k+\ell}$ on the (fixed) Cayley graph of $F_k$ that are associated with
the homomorphisms $\{h_n\}$,
converge after rescaling to
a faithful action of a (restricted) limit group $L$ on a real tree $Y$.

With the convergent subsequence of homomorphisms, for which the sequence of actions converges into the limit group $L$, we want to associate not just
an abelian  (JSJ like) decomposition, but also a $resolution$. i.e.,
a finite descending chain of (strict and proper) epimorphisms of limit groups, such that the terminal limit group of the resolution is a free
product of the coefficient group $F_k$ with a possibly trivial free group.

To get such a resolution, and afterwards use these resolutions to construct a Makanin-Razborov diagram,
we use the techniques that
appear in [Ja-Se] to construct a Makanin-Razborov diagram over free products.

The action of $L$ on the limit tree $Y$ can be analyzed using [Gu].
The action is superstable in the sense of [Gu],
stabilizers of tripods in the limit tree $Y$ are trivial, and stabilizers of non-degenerate segments are f.g.\ free abelian).
This enables us to apply the results of [Gu]
(theorem 4.1 in [Gu]), and deduce that either $L$ splits as a free product,
or the action of $L$ on $Y$ decomposes into a graph of actions
(see theorem 4.1 in [Gu] for these notions and conclusion).

If from the action of $L$ on $Y$ we obtain a free decomposition of $L$, we continue with each of the factors in that free decomposition in parallel.
We denote the abelian decomposition that is obtained from the action of $L$ on $Y$, $\Delta_L$. $\Delta_L$ contains  QH and (free) abelian vertex groups,
together with point stabilizers that are vertex groups,
and f.g.\ free abelian (possibly trivial) edge groups.

We continue by associating a modular group of automorphisms of $L$ with the decomposition $\Delta_L$, as it appears in section 5 in [Se1].
We denote the shortened homomorphisms  $sh_n=h_n \circ \varphi_n$, where $\varphi_n$ is a shortening automorphism from the modular group of $\Delta_L$.
If the sequence
$\{sh_n\}$ is a divergent sequence, then a subsequence
converges to an action of a limit group $L_2$ on a real tree $Y_2$, where $L_2$ is a quotient
(but not necessarily a proper quotient) of the limit group $L$.

Let $\Delta_{L_2}$ be the virtually abelian decomposition that is associated with the action of $L_2$ on $Y_2$, precisely as we associated
the virtually abelian decomposition $\Delta_L$ with the limit group $L$.
If $L_2$ is a proper quotient of $L$, we continue to the next step, precisely as we did in the construction of the action of $L_2$ on $Y_2$, starting with the
action of $L$ on $Y$. Since the action of $L_2$ on $Y_2$ was obtained from shortening the homomorphisms $\{h_n\}$ using the modular automorphisms
that are associated with $\Delta_L$, if
$L$ is isomorphic to $L_2$, $\Delta_{L_2}$ contains virtually abelian decompositions that do not appear in $\Delta_L$. i.e., the virtually abelian
decomposition in $\Delta_{L_2}$ give decompositions of rigid vertex groups in $\Delta_L$ that are compatible with $\Delta_L$, or it gives virtually abelian
decompositions that are hyperbolic w.r.t. $\Delta_L$.

In both cases it is possible to construct a JSJ type refinement of both $\Delta_L$ and $\Delta_{L_2}$, $\Delta^2_{L_2}$, from which it is possible to
extract both decompositions. Since $\Delta^2_{L_2}$ is a proper refinement of $\Delta_L$, either the free decomposition
of $L$ has changed, or the modular group that is associated with $\Delta^2_{L_2}$  is strictly bigger than
the modular group that is associated with $\Delta_L$. Hence, we can now repeat the construction of the action of $L_2$ on $Y_2$, by shortening with
automorphisms from the bigger modular group that is associated with $\Delta^2_{L_2}$.

As long as the shortening quotients are isomorphic to the limit group $L$, we get a sequence of proper refinements of the associated virtually
abelian decompositions. By Bestvina-Feighn accessibility or alternatively by acylindrical accessibility this refinement procedure terminates after
finitely many steps.

Since a refinement process terminates after finitely many steps,
the iterative procedure replaces the limit group $L$ with a proper quotient after finitely many steps. Hence, by the d.c.c.\ for limit groups
(section 5 in [Se1]), the iterative procedure of shortenings follow by refinements terminates after finitely many steps, and when it terminates we are left
with a limit groups which is a free product of the coefficient limit group $F_k$ and (possibly) a free group.

Note that some of the epimorphisms along the constructed finite sequence  are isomorphisms, and some are proper quotients.
The modular groups that are
associated with the limit groups  along the resolution are modular groups that are defined using the gradually refined  abelian decomposition of the
corresponding limit group, are dominated by  the modular groups that are associated with the JSJ decomposition of the corresponding limit group by
they may be strictly dominated by them.

As in [Se1] we call the constructed finite sequence of epimorphisms of limit groups a $resolution$. By construction the
constructed resolution has to be strict (definition 5.9 in [Se1]). The strictness is crucial for the rest of construction of the canonical MR diagram.

\vglue 1.5pc
\proclaim{Theorem 1.1} Let $U$ be a subset of $(F_k)^{\ell}$ and let $\{h:F_{k+\ell} \to F_k\}$ be the set of homomorphisms that are associated with $U$.
From any sequence of  homomorphisms that are associated with $U$ it is possible to extract a subsequence that  factors through
a resolution:
$L_1 \to L_2 \to \ldots \to L_f$.
\endproclaim

Note that if no refinements are applied along the levels of the constructed resolution then the sequence of homomorphisms and their
shortenings actually converge into the limit groups along the resolution and their associated abelian decompositions. Hence, such convergence
occurs when all the maps along the resolution are proper quotients between its limit groups, as then no refinements are needed.

In [Se2] we have introduced the $completion$ of a resolution (definition 1.12 in [Se2]). Note that by construction, the original limit group $L$,
which is the limit of the sequence of homomorphisms $\{h_n\}$,
embeds in the completion of the resolution $Res$: $\rho:L \to Comp(Res)$.

In a similar way to  what is proved in [Ja-Se]
for homomorphisms into free products, our next goal is to prove that  given a set $U$ there exist finitely resolutions, such that every homomorphism
$h: F_{k+\ell} \to F_k$ that is associated with an $\ell$-tuple from $U$, $h$ factors through one of the finitely many resolutions.

Given a set $U$
we start with all the possible sequences of homomorphisms, $h_n: F_{k+\ell} \to F_k$, that are associated with $\ell$-tuples from $U$, all their convergent
subsequences, and the resolutions that are constructed from such a convergent subsequence. Since each resolution is finite and limit groups are f.p.\ there
are at most countably many such resolutions.

We order the countable set of  resolutions that were constructed from convergent sequences of homomorphisms that are associated
with $\ell$-tuples in $U$.
If there are homomorphisms that do not factor through the  first $r$ resolutions,  we choose
a homomorphism that does not factor through the first $r$ resolutions.

If no finite subset of resolutions
suffice, we get an infinite sequence of homomorphisms, $\{h_r\}$, such that for each $r$ $h_r$ does not factor through the first $r$ resolutions.
Given this sequence, we can pass to a subsequence that converges into a resolutions.

This  resolution
 must appear in our ordered list of resolutions. Suppose that its place in the ordered list of resolutions is $r_0$.
Then for large indices $r>r_0$, homomorphisms $h_r$ from the convergent subsequence
do factor through this fixed resolution.  But this contradicts the choice of the homomorphisms
$\{h_r\}$, since they were supposed not to factor through resolutions that appear in the ordered list in the first $r$ places.

\vglue 1.5pc
\proclaim{Theorem 1.2} Let $\ell$ be a positive integer and $U$ be a set of $\ell$-tuples. There exists finitely many resolutions such that:
\roster
\item"{(1)}" all the homomorphisms $h:F_{k+\ell} \to F_k$ that are associated with $\ell$-tuples from $U$ factor through at least one of the finitely many
resolutions.

\item"{(2)}" for each of the finitely many resolutions, there exists a sequence of homomorphisms that are associated with points in $U$ that converge into
the limit groups along the resolutions and to their associated abelian decompositions (after finitely many  refinings).

\item"{(3)}" with each of the finitely many resolutions there is an associated map: $\rho_i: L_i \to Comp(Res_i)$, where $L_i$ is the limit of the original
sequence of homomorphisms $\{h_n\}$ that were used to construct the resolution.
\endroster
\endproclaim

\medskip
Theorem 1.2 associates a finite collection of resolutions with a set of $\ell$-tuples from a free group. However, the finite collection is not canonical.
To obtain a canonical collection we define a complexity of a finite collection of resolutions that guarantees that there exists a minimal
complexity collection, and prove that the minimal complexity collection is essentially unique.

In this paper we will analyze only sets of tuples that have a $single$ $ended$ MR diagrams. A single ended MR diagram is an MR diagram
that contains no free decompositions along its
levels, and for which no s.c.c.\ in a  QH vertex group along its levels is mapped to a trivial element in the next level.
Such (single ended) sets of tuples are important in studying varieties over free semigroups and associative algebras.
General sets will be analyzed in the second paper in the sequence. To define the complexity of a single ended diagram we first need to associate
a $modeled$ $structure$ with a resolution.

\proclaim{Definition 1.3} Let $Res$ be a  strict single ended   resolution
Let $Comp(Res)$ be its completion (definition 1.12 in [Se2]).
With the completion $Comp(Res)$ we associate a $modeled$ structure.

The modeled structure of the completion is composed from a collection of subtowers (subcompletions) that are part of the completion. Each subtower is
built over a QH vertex group along the completion, which is the  bottom (base) level of the subtower. Over the bottom QH vertex groups the subtower may contain
QH and
abelian vertex groups. Each QH vertex group is mapped isomorphically onto a suborbifold of the bottom QH vertex group. Each (free)  abelian vertex group
is of rank 2, where the  cyclic
edge group that is connected to it is mapped onto a s.c.c.\ in the bottom QH vertex group, and the  abelian vertex group
is generated by its  cyclic edge group and a formal Dehn twist generator that is added in the construction of the completion  (see section 1 in [Se2]).

The modular groups of all the QH vertex groups and all the  abelian vertex groups in a subtower are contained in the modular group of the
QH vertex group in the base of the subtower. Hence,
it is further required that when we collapse a subtower to the QH vertex group in its basis, using the retractions of the subtower,  the specializations
of the limit groups that is mapped into the (original) ambient tower and extend to specializations of the original ambient tower extend to specializations
of the collapsed tower.

We denote a modeled completion, $ModComp(Res)$. It consists of the completion and the form of the subtowers that are built over QH vertex groups along it.
Note that a subtower in the modeled structure can consist of a single QH vertex group, and not every  abelian vertex group in the completion is necessarily
contained in a modeled subtower.

We define the $reduced$ $modeled$ $completion$ to the completion that is obtained by collapsing each subtower in a modeled completion (using the retractions
in the modeled subtower) to
its bottom QH vertex group. We denote the reduced completion $Red(Comp(Res))$.

Note that there is an embedding of a limit group $L$, that is associated with the limit of the original sequence of homomorphisms $\{h_n\}$,
 into the completion of the original resolution, $Comp(Res)$, and this embedding restricts to an embedding of $L$ into the completion of a modeled resolution
and into its reduced modeled completion.
\endproclaim

\proclaim{Definition 1.4}
The complexity of a modeled (single ended) completion, $ModComp(CRes)$, is a tuple of integers ordered lexicographically in the following order:

\roster
\item"{(1)}" with each QH vertex group $Q$ that appears in the bottom of a modeled QH vertex group in one of the levels of the completion, $ModComp(CRes)$,
we associate a pair ($-\chi(Q)$,$g(Q)$), the Euler characteristic and the genus of the associated orbifold. The highest term in the complexity
of $ModComp(CRes)$ is the list of pairs that are associated with
these  QH vertex  groups along $ModComp(CRes)$, in a non-increasing lexicographical order.


\item"{(2)}" The sum $\Sigma \, rk(A_i) - rk(E_i)$ of the ranks of the  abelian vertex groups along the modeled completion minus the ranks of the
edge groups that are connected to them, where the sum is taken only over  abelian vertex groups that are not part of a subtower over a QH
vertex group in the modeled completion.

\item"{(3)}" The number of  abelian vertex groups along the completion (only those that are not part of a subtower over a QH vertex group).
\endroster
\endproclaim

Note that the complexity of a modeled completion is equivalent to the complexity of its reduced completion. Also,
note that complexities of modeled completions are well-ordered. To associate a canonical MR diagram with a set of tuples we define the
complexity of a (modeled) MR diagram,  study minimal
complexity MR diagrams, and prove that such MR diagrams are unique under a natural equivalence relation.

Complexities of single ended modeled completions extend to complexities of single ended (modeled) MR diagrams, i.e.,
to diagrams in which all the (modeled) resolutions are single
ended.

\proclaim{Definition 1.5} Let $\ell>0$ be an integer, and let $U$ be a set of $\ell$-tuples. Suppose that a single ended modeled MR diagram can be associated with
$U$  (cf. theorem 1.2). We define the complexity of the MR diagram to be the complexities of the resolutions in it ordered in a non-increasing order. Since complexities
of resolutions are well-ordered, the lexicographical order on the set of complexities of MR diagrams is well-ordered as well.
\endproclaim

\medskip
Theorem 1.2 associates a (non-unique) MR diagram with a set of tuples. Definitions 1.3 and 1.4 associate complexities with modeled resolutions, and not with ordinary resolutions.
On an ordinary resolution one can define a trivial modeled structure, by setting every QH subgroup along the resolution to be the modeled subtower over itself.
However, there is a natural procedure that associates with each resolution, or completion, an associated (most economical) modeled structure.
This procedure is needed in studying minimal complexity (modeled) MR diagrams.

Let $Res$ be a resolution from an MR diagram.
We go along the levels of the completion, $Comp(Res)$, from bottom to top,
and push down
QH and
abelian vertex groups, for which the  abelian edge groups that are connected to them  are all elliptic in the lower level of the completion.
We also push down
abelian vertex groups $A$ that are connected to edges with  cyclic edge groups $C$, such that $A$  is generated by $C$ and an additional
element that is associated with Dehn twists along $C$,  and QH vertex groups $Q$, that are mapped isomorphically into a s.c.c.\ or a subsurface of
a QH vertex
group in a lower level $\hat Q$.

To push down such a QH vertex group $Q$ or an  abelian vertex group $A$ we also require the following. Let $\Lambda_{\hat Q}$ be the
abelian decomposition that the QH vertex group $\hat Q$ inherits from the embedding of $Q$ into it, or from the splitting of it along the s.c.c.\ $C$,
which is the peg of an  abelian vertex group $A$ in a higher level.
Let $M_{\hat Q}$ be the limit group that is associated with all the levels in the completion, $Comp(Res)$ up to the level of $\hat Q$. $M_{\hat Q}$ has
a tower structure that it inherits from $Comp(CRes)$. With $M_{\hat QH}$ we associate an abelian decomposition, $\Lambda_{M_{\hat Q}}$,
that is obtained from the tower structure
of $M_{\hat Q}$ by collapsing all the QH and  abelian vertex groups in the tower structure that is associated with $M_{\hat Q}$ (and
inherited
from $Comp(Res)$) except $\hat Q$, and further replacing $\hat Q$ with
its abelian decomposition $\Lambda_{\hat Q}$, and collapse all the vertex groups and edge groups in $\Lambda_{\hat Q}$ except the ones that are
associated with the isomorphic image of the QH vertex group $Q$ or the edge group (s.c.c.\ ) that is connected to the  abelian vertex group $A$.

We require that all the  abelian and the QH vertex groups  in the completion, $Comp(CRes)$,  that are in levels above the level of $\hat Q$ and up to the
level of $A$ or $Q$ in $Comp(Res)$ are mapped to elliptic subgroups in $\Lambda_{M_{\hat Q}}$, by the appropriate compositions of the retraction maps
between consecutive levels in the completion, $Comp(Res)$.

If there are edge groups or QH vertex groups in some level of a cover resolution, that are mapped isomorphically into a QH vertex group in a lower level,
and satisfy these conditions, we push them down to be part of a tower over the QH vertex group $\hat Q$ in a lower level, which is going to be part of the modeled
structure of the obtained cover resolution.

After pushing down all the  abelian and QH vertex groups that satisfy the above conditions for a push down, we continue to the next upper level.
In analyzing the next level, we replace the QH vertex group $\hat QH$ that was part of some lower level of the completion, with the subtower that is associated
with it, in which it is its bottom level.

We continue climbing up along the levels of the completion, $Comp(Res)$, by pushing down  abelian vertex groups and QH vertex groups
if all of their edge groups are elliptic in a level below them or if they satisfy the conditions to
be part a tower over a QH vertex group in a lower level, taking into account that QH vertex groups in lower levels have changed according to the modeled structures
they already have (i.e., according to the subtowers that are already associated with them from pushing down QH and virtually abelian vertex groups from higher
levels).

After climbing to the top level, and if along climbing from the bottom level
there were QH or abelian vertex groups that were pushed to a lower level or pushed to be part
of a (modeled) subtower
above a QH vertex group in a lower level, we continue the procedure by climbing again from the bottom level and taking into account the new
abelian decompositions and the new modeled subtowers.

This  procedure (of climbing and pushing down QH and abelian vertex groups)  terminates after finitely many steps,
since the number of edges and vertices in the various levels
of the original completion, $Comp(Res)$,  is finite, and each vertex group can move to a lower level only boundedly many times (i.e., at most the number of levels
in the completion).

After this procedure
terminates, all the edge groups in the modified completion are virtually cyclic, and the terminal completions are is modeled.
An edge group that is not part of a subtower over a QH vertex group in the  modeled completion, and that is not connected to a QH
vertex group in an  abelian decomposition that is associated with a level above
the two terminating ones in the terminal completion, is not elliptic in the  abelian decomposition that is associated with the level below it in  the
terminal completion.
Furthermore,
at least one of  the  cyclic
edge groups that are connected to each of the  QH vertex groups, that appear in the bottom of each subtower over a QH vertex group in the modeled
structure of the terminal completion,  and the QH vertex group is in an  abelian decomposition that is associated with a level above the two
terminating ones in $Res$,
is not elliptic in the the abelian decomposition that is associated with the  level below the level in which the QH vertex group appears in .the terminal
completion.

To analyze minimal complexity resolutions and diagrams we will need the constructed resolutions to have weak test sequences from the given tuples. Test sequences
were defined in [Se2] and are used there for analyzing sentences and formulas with 2 quantifiers. Test sequences satisfy strong combinatorial properties
and in general we can not expect that resolutions that are associated with a set of tuples (as in theorem 1.2) will have test sequence of homomorphisms that
associated with the given set of tuples. Hence, we need to considerably relax the properties of the sequences that we are going to work with.

\vglue 1.5pc
\proclaim{Definition 1.6} Let $h_n: L \to F_k$ be a sequence of homomorphisms from a limit group $L$ to the free group $F_k$, that factor
through a resolution  $Res$ of $L$. The homomorphisms $h_n$ naturally extend to homomorphisms $\tilde h_n: Comp(Res) \to F_k$.

We say that the sequence of homomorphisms $\{h_n\}$ is a $weak$ $test$ $sequence$ of the resolution $Res$, if the sequence $\{\tilde h_n\}$ converges into
a faithful action of $Comp(Res)$ on a real tree, from which it is possible to extract the abelian decompositions along the various levels of the completion
$Comp(Res)$. i.e., this means that every QH vertex group along the completion acts on an IET component of a corresponding real tree, and each abelian
vertex group acts discretely or axially on a component of such a tree.
\endproclaim

If a resolution has a weak test sequence of homomorphisms that factor through it, then the same sequence of homomorphisms factors through the
modeled completion that is obtained from the original resolution using the iterative procedure that we presented, and the original sequence is
a weak test sequence of the obtained modeled completion.

The construction of a resolution in theorem 1.1 guarantees that if all the maps along the constructed resolution are proper quotients then the sequence
of homomorphisms from which the resolution was constructed is a weak test sequence. Hence, in this case the sequence of homomorphisms is also a weak test
sequence of the modeled completion that is constructed from such a resolution.

However, if a resolution that is constructed in theorem 1.1 contains isomorphisms (and not only proper quotients),
then the sequence of homomorphisms from which the resolution was constructed need not be a weak test sequence. When we analyze set of tuples and their associated
MR diagrams we need that all the resolutions in the MR diagram will have weak test sequences of homomorphisms that are associated with tuples from the
given set.

Since the construction of resolutions in theorem 1.1 does not guarantee that the sequence from which a resolution is constructed is a weak test sequence
we need a variation of the construction  that will guarantee this property. Indeed such a variation exists in our work on word equations [Se9], i.e., in
our work on Makanin-Razborov diagrams for a system of equations over a free semigroup.

The Makanin-Razborov diagram that is constructed in [Se9] is for semigroups and not for groups, but it works easier for set of tuples in groups,
indeed without the
complications of positivity that exist in analyzing equations over semigroups. The resolutions that are constructed there are finite, but QH vertex groups
along them may be obtained from resolutions with countably many levels (in fact with each QH vertex group it is possible to associate a countable
ordinal that is associated with a modular tower above it).

By the construction of resolutions that is described in [Se9], the sequence of
homomorphisms from which the countable resolutions are constructed are weak test sequences for these countable resolutions. The complexity
of such a countable resolution is the complexity of a finite (reduced) resolution that is obtained by collapsing the countable subtowers
to the QH vertex groups at their bases. Therefore, for the rest
of this section we will continue to work with these countable resolutions, and  MR diagrams that contain finitely many such resolutions that are
constructed from homomorphisms that are associated with tuples from a given set of tuples, and such that every homomorphism that is associated with such a tuple
factors through at least one of the resolutions from the MR diagram (cf. theorem 1.2).

\vglue 1.5pc
\proclaim{Theorem 1.7} Let $\ell$ be a positive integer and let $U$ be a set of $\ell$-tuples. There exists finitely many (countable)
modeled resolutions such that:
\roster
\item"{(1)}" each homomorphism $h:F_{k+\ell} \to F_k$ that is associated with an $\ell$-tuple from $U$ extends to a homomorphism that
factors through the reduced completion of
at least one of the (finitely many)
resolutions  of the (countable) modeled resolutions.

Note that a homomorphism, $h:F_{k+\ell} \to F_k$, factors through the reduced completion of a single ended resolution $Res$, if it extends to a homomorphism of the
reduced completion, $\hat h:Red(Comp(Res)) \to F_k$, and no s,c,c,\ in any of the QH vertex groups along the reduced completion is mapped by $\hat h$
to the identity in $F_k$.

\item"{(2)}" for each of the finitely many (countable) modeled resolutions there exists a sequence of homomorphisms that are associated with tuples in $U$ that form a weak
test sequence for (the completion of) the (countable) modeled resolution.
\endroster
\endproclaim


%

In this paper we will analyze only single ended diagrams. i.e., we will assume that all the resolutions in a collection that satisfies
the conclusions of theorem 1.7 w.r.t.\ a set of tuples $U$
are single ended. Note that by the construction in [Se9] if a set of tuples $U$ has a finite collection of single ended resolutions that satisfy the properties of
theorem 1.2, then there is such a finite collection of single ended (countable) modeled resolutions that satisfies the conclusions of theorems 1.7.

Let $U$ be a set of tuples and suppose that it has a finite collection of single ended modeled resolutions that satisfies the conclusion of theorem 1.7. We look
at all the finite single ended collections that satisfy the conclusions of theorem 1.7 w.r.t.\ the set $U$. Since complexities of resolutions and diagrams are
well-ordered there exist finite (single ended) collections that satisfy the conclusions of theorem 1.7 with minimal complexity.
Our goal is to prove that minimal complexity (single ended)
collections are essentially
unique. i.e., that minimal complexity (single ended) collections belong to the same natural equivalence class of collections.

For presentation purposes, we start by analyzing  minimal complexity diagrams in case the abelian decompositions along all the levels of  the
completions of the modeled resolutions
in a minimal complexity
(single ended) diagram of a set of tuples $U$, contain no abelian vertex groups that are not part of a modeled subtower above some QH vertex  group in the modeled resolution
(in the case of existence of  abelian groups that are not part of modeled subtowers, to obtain uniqueness of the diagram we will need to slightly modify
part (1) from the properties that an MR diagram needs to satisfy in the statement of theorem 1.7 (see definition 1.15 below)).

\vglue 1.5pc
\proclaim{Theorem 1.8} Let $\ell$ be a positive integer and let $U$ be a set of $\ell$-tuples in $F_k$. Suppose that $U$ has two minimal complexity
single ended finite collections of
(countable) modeled resolutions that satisfy the conclusions of theorem 1.7, and suppose that the abelian decompositions along the levels of the modeled resolutions in these two minimal
complexity diagrams contain no abelian vertex group that is not part of a modeled subtower above some QH vertex group.

Let $Res_1,\ldots,Res_m$ be the modeled resolutions in the first collection, and $Res^1_1,\ldots,Res^1_t$ be the modeled resolutions in the second collection.
Recall that with each resolution there is an associated embedding of a limit group $L_i$ into the reduced completion of the resolution $Res_i$, and a
similar associated embedding of a limit group $L^1_i$ into the
reduced completion of $Res^1_i$.

Then $m=t$, and up to a change of order,for each $i$, $1 \leq i \leq m$, the structures of the reduced
completions, $Red(Comp(Res_i))$ and $Red(Comp(Res^1_i))$, are similar (same number of levels
and same number and isomorphic  QH  vertex groups  in each level), and there exist maps:
$$\eta_i:Red(Comp(Res_i)) \ \to \ Red(Comp(Res^1_i)) \ \ and \ \
\eta^1_i:Red(Comp(Res^1_i)) \to Red(Comp(Res_i))$$
that map the QH vertex groups in one completion isomorphically onto  corresponding QH vertex groups in the target completion, and the image of the limit group $L_i$
that embeds in $Red(Comp(Res_i))$ isomorphically onto the image of $L^1_i$ in $Red(Comp(Res^1_i))$ and vice versa.

Hence, the pairs of maps, $\eta_i$ and $\eta^1_i$, are isomorphisms that map $QH$ vertex groups isomorphically onto QH vertex groups   in the completions
of the completions of the corresponding reduced modeled
completions.
\endproclaim

\nfp Given a resolution from the diagram we start by trying to push down  QH vertex groups using the iterative procedure that we presented. After
finitely many steps this procedure terminates.
After the procedure  terminates,  in a QH vertex group $Q$ that is not from the bottom two levels of
a modeled resolution from a minimal complexity diagram, and is in the bottom of a modeled subtower,
at least one of its boundary subgroups is not elliptic in the abelian decomposition that
is associated with the level below it, and $Q$ can not be pushed to be part of a model subtower that has its bottom below the level of $Q$ in its resolution.

For presentation purposes suppose first that the modeled structures of the modeled resolutions in the two minimal complexity MR diagrams are trivial, i.e.,
that the modeled subtowers above any QH subgroup in these modeled resolutions contain only the base QH vertex groups with nothing above it. In this case
the reduced modeled resolutions are identical to the modeled resolutions.

Let $Res_1$ be a resolution of maximal complexity from the resolutions, $Res_1,\ldots,Res_m$. By the properties of the MR diagrams that are associated
with the set $U$, there exists a
sequence of homomorphisms $\{h_n\}$ that are associated with elements of the set $U$ and extend to a weak test sequence of $Res_1$.

Each of the homomorphisms $\{h_n\}$ extend to homomorphisms that factor through one of the resolutions, $Res^1_1,\ldots,Res^1_t$. By the construction of
formal solutions and formal limit groups in [Se2],
from the sequence $\{h_n\}$ and its extensions to $Comp(Res_1)$ and to the completions of the resolutions, $Res^1_1,\ldots,Res^1_t$ it is possible
to pass to a further subsequence that converges to a map:  $\nu:Comp(Res^1_i) \to FCl(Res_1)$, where $FCl(Res_1)$ is a $framed$  closure of $Comp(Res_1)$,
and $\nu$ maps the image of the limit group $L^1_i$ in $Comp(Res^1_i)$ onto the image of the limit group $L$ in $FCl(Res_1)$.

Framed closures were introduced in definition 6 in [Se5].
A framed closure of $Comp(Res_1)$ (that does not contain abelian vertex groups) is obtained from the completion, $Comp(Res_1)$,
by possibly replacing some of the QH vertex groups along it with supergroups
of finite index, i.e., replace the surfaces that are associated with these QH vertex groups by surfaces that are finitely covered by them.

\proclaim{Lemma 1.9} With the notation above if $FCl(Res_1)$ is not identical to $Comp(Res_1)$ then either:
\roster

\item"{(1)}" two non-conjugate  QH vertex groups in $Comp(Res_1)$ are mapped into conjugate ones
in $FCl(Res_1)$.

\item"{(2)}" Some QH vertex groups in the framed closure $FCl(Res_1)$ properly contain QH vertex groups in  $Comp(Res_1)$.
\endroster
\endproclaim

\nfp Follows from the construction of a closure in [Se2].

\line{\hss$\qed$}

Note that in both cases in lemma 1.9 the complexity of $FCl(Res_1)$ is strictly smaller than the complexity of $Comp(Res_1)$.
Now, suppose that
$FCl(Res_1)$ is identical to $Comp(Res_1)$.

\proclaim{Lemma 1.10} If $FCl(Res_1)$ is identical to $Comp(Res_1)$ then either:
\roster
\item"{(1)}" The map $\nu:Comp(Res^1_i) \to FCl(Res_1)$ is an isomorphism. In particular, the two completions $Comp(Res_1)$ and $Comp(Res^1_i)$ have the same
structure and the homomorphisms $\{h\}$ that are associated with tuples from $U$ and extend to homomorphisms that factor through $Comp(Res_1)$ extend to
homomorphisms that factor through $Comp(Res^1_i)$ and vice versa.

\item"{(2)}"  The map $\nu$ is not onto. In that case there exists a QH vertex group $Q$ in $FCl(Res_1)$ that the image $\nu(Comp(Res^1_i))$ intersect
trivially all its conjugates. i.e., it intersects all its conjugates only in conjugates of some boundary subgroups.
\endroster
\endproclaim

\nfp
The argument that we use to prove lemma 1.10 is based on the analysis of quotient resolutions of minimal rank, that appears in section 1 of [Se3],
and in section 1 in [Se4].

Let $Q_1$ be a QH vertex group in the terminal  abelian decomposition  of $Comp(Res^1_i)$. All the boundary components of $Q_1$ can be conjugated
into the terminal level of $Comp(Res^1_i)$.
By construction, $\nu$ maps the boundary elements of $Q_1$ into the terminal  subgroup in $FCl(Res_1)$. Suppose that $\nu(Q_1)$ intersects
conjugates of some QH vertex groups  along the levels of $FCl(Res_1)$ non-trivially. i.e.,
not only in conjugates
of boundary subgroups. Suppose that  $Q_u$ is a QH vertex group in $FCl(Res_1)$, that intersects non-trivially a conjugate of $\nu(Q_1)$,
and $Q_u$ is in the  highest level in $FCl(Res_1)$ for which there are such $QH$  vertex groups.

Since we assumed that the MR diagrams are single ended, a subsurface of $Q_1$ is mapped by $\nu$ onto a finite index subgroup of $Q_u$.
Therefore,
the complexity (pair) of $Q_u$ is bounded by the complexity of $Q_1$, with equality if and only if  $Q_u$ is a QH vertex group in the terminal  abelian
decomposition of $FCl(Res_1)$, and $\nu$ maps $Q_1$ isomorphically onto a conjugate of $Q_u$.

Let $j$ be the highest level in  $FCl(Res_1)$ for which $\nu(Q_1)$ intersects non-trivially a conjugate of a
QH vertex group.
Let $FCl(Res_1)_{j+1}$ be the subgroup which is the part of $FCl(Res_1)$
that contains all the levels up to level $j+1$ (the level below level $j$).  Let $\eta_{j+1}$ be the retraction: $\eta_{j+1}:FCl(Res_1) \to FCl(Res_1)_{j+1}$.

We look at the image of $Q_1$ under the composition: $\eta_{j+1} \circ \nu$. The conjugacy classes of the
images of boundary components of $Q_1$ do not change by the composition with
the retraction $\eta_{j+1}$. Suppose that the image of $Q_1$ under the composition intersects non-trivially
a conjugate of a QH vertex group $Q^{j+1}$ in level $j+1$
of $FCl(Res_1)$, which is
the highest level in the part of $FCl(Res_1)$ that is the image of the retraction $\eta_{j+1}$.

By the same arguments that applied to the QH vertex group
$Q_u$, the complexity of $Q^{j+1}$ is bounded by the complexity of $Q_1$, with equality
if and only if  $Q^{j+1}$ is a QH vertex group in the terminal  abelian
decomposition of $FCl(Res_1)$, and $\eta_{j+1} \circ  \nu$ maps $Q_1$ isomorphically onto a conjugate of $Q^{j+1}$.

We continue iteratively, by composing with retractions to lower levels. By the same arguments we conclude that any QH vertex group that has a conjugate
that is intersected non-trivially
by the image of $Q_1$ under the composition of a retraction and the map $\nu$, has complexity that is bounded by the complexity of $Q_1$, with equality
if and only if the QH vertex group is in the bottom level of $FCl(Res_1)$, and the composition of the retraction and $\nu$ maps $Q_1$ isomorphically
onto that QH vertex group.

\smallskip
The bounds on the complexities of the QH vertex groups in $FCl(Res_1)$ that intersect non-trivially conjugates of the images of a
QH vertex group $Q_1$ in the bottom  abelian decomposition of $Comp(Res^1_i)$, under compositions of retractions of $FCl(CRes)$ with the map
$\nu$, enable us to analyze the image of $\nu$.

If the complexities of all the QH vertex groups in $FCl(Res_1)$ that have conjugates that are intersected non-trivially by $\nu(Q_1)$ and its
images under the retractions,
are strictly smaller  than the complexity of $Q_1$,
we gain in the ambient
complexity. i.e.,  $\nu(Q_1)$ covers only QH vertex groups in $FCl(Res_1)$ that have strictly smaller complexity.

If there is a QH vertex group $Q$ in $FCl(Res_1)$ that has a conjugate that is intersected non-trivially by $\nu(Q_1)$ or its images under retractions,
and $Q$
has the same complexity as $Q_1$, $Q$  has to be in the terminal  abelian decomposition of $FCl(Res_1)$, and the composition of the map $\nu$
with (possibly) retractions of $FCl(Res_1)$ have to map $Q_1$  onto a conjugate
of $Q$ in $FCl(Res_1)$.

If there are two QH vertex groups $Q_1$ and $Q_2$
in the terminal  abelian decomposition of $Comp(Res^1_i)$, such that the images of them under a composition of a retraction of $Fcl(Res_1)$
with $\nu$ intersect conjugates
of the same QH vertex group $Q$ in the terminal  abelian decomposition of $FCl(Res_1)$, we also gain in the ambient complexity.

\proclaim{Definition 1.11 (cf. definition 1.6 in [Se3])} A QH vertex group $Q$ in the terminal virtually abelian decomposition of $FCl(Res_1)$
is called $surviving$ $surface$ if:

\roster
\item"{(i)}"  there exists
a QH vertex group $Q_1$ in the terminal  abelian decomposition of $Comp(Res_1)$ such that
the composition of a retraction with $\nu$ maps $Q_1$ isomorphically onto a conjugate of $Q$.

\item"{(ii)}" the images of all the other QH vertex groups in the terminal level of $FCl(Res_1)$ under the composition of a retraction
of $Fcl(Res_1)$ and $\nu$,
intersect every conjugate of $Q$ trivially (i.e., in a coin jugate of a boundary subgroup)
\endroster
\endproclaim

Suppose that $Q$ is a surviving surface, and let $Q_1$ be the unique QH vertex group in the terminal  abelian decomposition
of $Comp(Res^1_i)$, that is mapped by a composition of (possibly) a retraction of $Fcl(Res_1)$ and $\nu$ isomorphically onto a conjugate of $Q$.
We continue in a similar way to what was done
in section 1 in [Se3] for bounding the complexity of quotient resolutions.

We  rearrange the QH vertex groups
along the levels of the  closure, $FCl(Res_1)$, that are mapped by retractions isomorphically onto subsurfaces of the surviving surface  $Q$,
and modify accordingly the map $\nu$. We do that to guarantee that for each surviving surface $Q$
in the terminal  abelian decomposition of $FCl(Res_1)$,
the unique QH vertex group $Q_1$ that was mapped isomorphically onto a conjugate of $Q$  by a composition of a retraction and the map $\nu$,
will be mapped isomorphically onto $Q$ by the modification of $\nu$ (with no composition with a retraction). This will imply that a QH vertex group
$Q_1$ in $Comp(Res^1_i)$ that is mapped onto a conjugate of a surviving surface $Q$, intersects trivially conjugates of all the other QH vertex groups
in the modified $FCl(Res_1)$. Hence, the contribution of $Q_1$ to the ambient complexity of $Comp(Res^1_i)$ will be the same as the contribution of the
surviving surface $Q$ to
the complexity of $FCl(Res_1)$.

Let $Q^j,Q^{j+1},\ldots,Q^{r-1}=Q$, be the sequence of (not necessarily connected, some of which possibly empty) QH vertex groups (that we consider as surfaces)
along the  levels $j$ to $r-1$ of $FCl(Res_1)$'
such that the image of $Q_1$ under a composition of the corresponding
retraction with $\nu$ intersects non-trivially a conjugate of them.

Since $Q_1$ is mapped isomorphically onto a conjugate of the surviving surface $Q$
by a composition of a retraction and $\nu$, subsurfaces of $Q_1$ are mapped isomorphically onto conjugates of the QH vertex groups in the sequence:
$Q^j,\ldots,Q^{r-1}$, and s.c.c.\ in $Q_1$ are mapped onto  conjugates of edge groups in  abelian decompositions
along $FCl(Res_1)$. Furthermore, the final retraction maps QH vertex groups in $Q^j,\ldots,Q^{r-2}$ isomorphically onto conjugates of suborbifolds of
$Q$, and  the boundary subgroups of these QH vertex groups  onto s.c.c.\ in conjugates of $Q$.
Since $Q$ is a surviving surface, the images of all the other  QH vertex groups (except $Q_1$) in the terminal level of $Comp(Res^1_i)$, under compositions of
appropriate retractions of $FCl(Res_1)$ with $\nu$, intersect trivially conjugates of all the QH and virtually abelian vertex groups in the sequence:
$Q^j,\ldots,Q^{r-1}$.

These observations allow us to modify the structure of the reduced framed closure $FCl(Res_1)$ and the map $\nu$ in a similar way to the
modification of quotient resolutions
(in the minimal rank case)
in the first sections of [Se3] and [Se4].

We modify the closure $FCl(Res_1)$ by inverting the order of the levels of the QH  vertex groups that appear in the sequence:
$Q^j,\ldots,Q^{r-2}$, leaving $Q^{r-1}=Q$ unchanged. i.e., we set it to be: $Q^{r-2},\ldots,Q^j,Q^{r-1}=Q$.
This order inversion is completed by changing accordingly the map $\nu$,
that after modification maps $Q_1$ isomorphically onto a conjugate of $Q$ (without composing it with a retraction). We also modify accordingly the retraction
maps between consecutive levels, so that they are compatible with the reverse order of the QH vertex groups $Q^j,\ldots,Q^{r-2}$.
For the detailed modification of the map $\nu$ see section 1 in [Se3].
We denote the modified map $\tilde \nu$, and the modified reduced framed closure, with its modified retractions, $\tilde {FCl}(Res_1)$.

Note that the images of the other QH  vertex groups in the terminal  abelian decomposition  in  $Comp(Res^1_i)$ are not effected by
the modification of $\nu$, as these images intersect the QH and virtually abelian vertex groups in the sequence $Q^j,\ldots,Q^{r-1}=Q$ trivially.
The images of QH vertex groups in higher level  abelian decompositions in $Comp(Res_1)$ may be effected by the modification
of the map $\nu$, and the same holds for the image of $L$, $\nu(L)$, in the modified closure $\tilde {FCl}(Res_1)$ (see section 1 in [Se3] for
these modifications).

We perform these modifications of the  closure $FCl(Res_1)$, and of the map $\nu$,
for all the surviving surfaces. Note that because of trivial intersections,
these modifications can be conducted in parallel.

\smallskip
We continue iteratively by climbing along the levels of $Comp(Res^1_i)$,
and possibly further modifying the structure of
$\tilde {FCl}(Res_1)$ and the map $\tilde \nu$.
After completing the modification at level $b$, $2 \leq b \leq r^1_b-1$, we set $M_b$ to be the subtower of the   framed closure, $\tilde {FCl}(Res_1)$,
that contains the QH
vertex groups in $\tilde {FCl}(Res_1)$ that intersect non-trivially images (under $\tilde \nu$)
of QH  vertex groups in the
levels $b$ and below in the  completion $Comp(Res^1_i)$.
By construction,  $M_b$  is a subtower that is closed under the retractions of
the ambient modified closure $\tilde {FCl}(Res_1)$.

Note that the images under $\tilde \nu$ of all the edge groups in the  abelian decomposition in level $b-1$ in $Comp(Res^1_i)$, i.e., all the
boundary subgroups of QH vertex groups in this  abelian decomposition, are contained in $M_b$. Hence, if a QH vertex group $Q_{b-1}$
in the  abelian decomposition that is associated with level $b-1$ in $Comp(Res^1_i)$, and $Q$ is a QH vertex group in the highest possible level in
$\tilde {FCl}(Res_1)$
that
is not in $M_b$ and intersects non-trivially a conjugate of $\tilde \nu(Q_{b-1})$, then the complexity of $Q$ is bounded by the
complexity of $Q_{b-1}$.  Furthermore, equality in the complexities of $Q$ and $Q_{b-1}$ occurs if and only if $Q$ is a QH vertex group in some level of
$\tilde {FCl}(Res_1)$, and its boundary subgroups can be conjugated into $M_b$.

We look at compositions of $\tilde \nu$ and retractions of the closure, $\tilde {FCl}(Res_1)$, precisely as we did in analyzing QH vertex groups
in the terminal
abelian decomposition of $Comp(Res^1_i)$.

The analysis is complete when we analyze the highest level in $Comp(Res^1_i)$. The two given MR diagrams were assumed to be of minimal
complexity, and the resolution $Res_1$ is assumed to be of maximal complexity among the resolutions in the first diagram. Hence,  if every QH vertex group
in $FCl(Res_1)$
has a conjugate that intersects non-trivially the image $\tilde \nu (Comp(Res^1_i)$, then $Res^1_i$ has to be of maximal complexity in the second
MR diagram, every QH vertex group in $Comp(Res^1_i)$ is mapped isomorphically onto a QH vertex group in $\tilde {FCl}(Res_1)$, and distinct QH vertex groups in
$Comp(Res^1_i)$ are mapped isomorphically onto non-conjugate QH vertex groups in $\tilde {FCL}(Res_1)$. Since, by construction,  $\tilde {FCl}(Res_1)$
is isomorphic to $FCl(Res_1)$, we get that in this case $Comp(Res^1_i)$ is isomorphic to $FCl(Res_1)$. This is case (1) in lemma 1.9.

If there exists a QH vertex group in $\tilde {FCl}(Res_1)$ that intersects trivially the conjugates of the images (under $\tilde \nu$)
of all the QH vertex groups in
$Comp(Res^1_i)$, then the  image of $\tilde \nu$ intersects trivially the conjugates of a QH vertex group in $\tilde{FCl}(Res_1)$. This is case (2) in
lemma 1.9. In that case the image of $\tilde \nu$ is contained in a proper subtower of $\tilde {FCl}(Res_1)$, i.e., in a subtower that is obtained from
$\tilde {FCl}(Res_1)$ by omitting from it at least one QH vertex group.

If there isn't such a QH vertex group in $\tilde {FCl}(Res_1)$ then the map $\tilde \nu$ is necessarily onto, and by the Hopf property for
limit groups, $\tilde \nu$ is an isomorphism and case (1) in the lemma holds.

\line{\hss$\qed$}

So far we proved lemmas 1.9 and 1.10 assuming that the modeled resolutions in the two given minimal complexity MR diagrams, that are given in the
assumptions of theorem 1.8, have trivial modeled
structure. i.e.' that every modeled subtower contains a single QH vertex groups with nothing above it. At this point we continue
with the proof of theorem 1.8 and omit this assumption, i.e., we allow the
modeled structure to be arbitrary.

We use the notation as in the statement of theorem 1.8. By our assumptions, in case of general minimal rank modeled resolutions that
contain no abelian vertex groups that are not parts of modeled subtowers, each resolutions from the two diagrams has a weak test sequence of
homomorphisms that are associated with tuples in the subset $U$.

We assume that $\{h_n\}$ is a weak test sequence of the resolution $Res_1$, so by our assumptions these homomorphisms extend
to homomorphisms that factor through the reduced completions:
$Red(Comp(Res^1_1),\ldots,Red(Comp(Res^1_t))$. By the construction of
formal solutions and formal limit groups in [Se2],
from the sequence $\{h_n\}$ and its extensions to $Comp(Res_1)$,
it is possible
to pass to a further subsequence that converges to a map:  $\hat {\nu}:Red(Comp(Res^1_i)) \to FCl(Res_1)$, where $FCl(Res_1)$ is a $framed$
closure of $Comp(Res_1)$. By further composing the map $\hat {\nu}$ with the retractions of the modeled subtowers in $Comp(Res_1)$, we get a map:
$\nu: Red(Comp(Res^1_i)) \to FCl(Red(Comp(Res_1)))$.

The existence of the map $\nu$ enable us to continue along the proofs of lemma 1.9 and 1.10 in the case of resolutions with a trivial modeled structure
and generalize them to resolutions with a general modeled structure. For completeness we generalize the statements of lemmas
1.9 and 1.10 in the case of resolutions with a general modeled structure.

\proclaim{Lemma 1.12} With the notation above if $FCl(Red(Comp(Res_1)))$ is not identical to $Red(Comp(Res_1))$ then either:
\roster

\item"{(1)}" two non-conjugate virtually  QH vertex groups in $Red(Comp(Res_1))$ are mapped to conjugate ones
in $FCl(Red(Comp(Res_1)))$.

\item"{(2)}" Some QH vertex groups in the framed closure $FCl(Red(Comp(Res_1)))$ properly contain QH vertex groups in  $Red(Comp(Res_1))$.
\endroster
\endproclaim

\proclaim{Lemma 1.13} If $FCl(Red(Comp(Res_1)))$ is identical to $Red(Comp(Res_1))$ then either:
\roster
\item"{(1)}" The map $\nu:Red(Comp(Res^1_i)) \to FCl(Red(Comp(Res_1)))$ is an isomorphism. In particular, the two completions $Red(Comp(Res_1))$ and
$Red(Comp(Res^1_i))$ have the same
structure and the homomorphisms $\{h\}$ that are associated with tuples from $U$ and extend to homomorphisms that factor through $Red(Comp(Res_1))$ extend to
homomorphisms that factor through $Red(Comp(Res^1_i))$ and vice versa.

\item"{(2)}"  The map $\nu$ is not onto. In that case there exists a QH vertex group $Q$ in $FCl(Red(Comp(Res_1)))$ that the image $\nu(Red(Comp(Res^1_i)))$
intersect
trivially all its conjugates. i.e., it intersects all its conjugates only in conjugates of some boundary subgroups.
\endroster
\endproclaim

Lemmas 1.12 and 1.13 enable us to complete the proof of the uniqueness of a minimal complexity MR diagram under the assumptions of
theorem 1.8. We start with the modeled resolutions in the first minimal complexity MR diagram $Res_1,\ldots,Res_m$, and one of its maximal complexity
resolutions $Res_1$. By our assumptions, there exists a sequence of homomorphisms $\{h_n:F_{k+\ell}\to F_k\}$  that are associated with
$\ell$-tuples from the given set $U$, that extend to homomorphisms that form a weak test sequence for $Res_1$.

The weak test sequence of homomorphisms $\{h_n\}$ extend to homomorphisms that factor through the reduced completions of resolutions in the second
MR diagram $Red(Comp(Res^1_1)),\ldots,Res(Comp(Res^1_t))$. Suppose that a subsequence of the homomorphisms $\{h_n\}$ extend to homomorphisms that factor
through a reduced completion, $Red(Comp(Res^1_i))$.

By possibly passing to a further subsequence, and by the techniques of [Se2],
it is possible to obtain a map: $\nu:Red(Comp(Res^1_i)) \to FCl(Res(Comp(Res_1)))$. Suppose that
$FCl(Red(Comp(Res_1)))$ is identical to $Red(Comp(Res_1))$ and that the map $\nu$ satisfies part (1) of lemma 1.13, i.e., that $\nu$ is an isomorphism.
In that case the reduced completions, $Red(Comp(Res_1))$ and $Red(Comp(Res^1_i))$,  are isomorphic. Hence, every homomorphism $h:F_{k+\ell} \to F_k$ that
is associated with a point in $U$, and extends to a homomorphism that factors through $Red(Comp(Res_1))$, extends to a homomorphism that factors through
$Red(Comp(Res^1_i))$, and vice versa, since the two reduced completions are isomorphic.

Hence, in this case of isomorphic reduced completions in the two minimal complexity diagrams, we can restrict to the subset of
homomorphisms $h$, that are associated with
tuples in $U$, and do no extend to homomorphisms that factor through $Red(Comp(Res_1)$ (or $Red(Comp(Res^1_i))$). We denote this collection of
homomorphisms $U_1$. Since the two given diagrams are of
minimal complexity, if we take out $Res_1$ from the first diagram and $Res^1_i$ from the second diagram, the obtained diagrams must be minimal
complexity diagrams for homomorphisms from $U_1$. Therefore, in this case we are able to continue iteratively.

Suppose that  $FCl(Red(Comp(Res_1))$ is not identical to $Red(Comp(Res_1))$, i.e., that part (1) or (2) in lemma 1.12 holds, or that the map
$\nu:Red(Comp(Res^1_i)) \to FCl(Red(Comp(Res_1)))$ is not onto. i.e., that part (2) in lemma 1.13 holds. In that case we look at all the
sequences of homomorphisms $\{h_n\}$ that are associated with tuples from the set $U$, that satisfy:
\roster
\item"{(1)}" the homomorphisms $\{h_n\}$ do not extend to homomorphisms that factor through the reduced completions of $Res_2,\ldots,Res_m$. Hence,
they must extend to homomorphisms that factor through the reduced completion of $Res_1$.

\item"{(2)}" the homomorphisms $\{h_n\}$ extend to homomorphisms
that form a weak test sequences of a resolution $CRes$.
\endroster

Using the techniques of [Se2] that we applied above, by possibly passing to a further subsequence, from the weak test sequence $\{h_n\}$ it
is possible to construct a map: $\tau:Red(Comp(Res_1)) \to FCl(Red(Comp(CRes)))$. Suppose that $FCl(Red(Comp(CRes)))$ is identical to $Red(Comp(CRes))$,
and that the map $\tau$ satisfies property (1) in lemma 1.13, i.e., suppose that $\tau$ is an isomorphism.

The weak test sequence of homomorphisms $\{h_n\}$ extend to homomorphisms that factor through the reduced completions of the resolutions in
the second given MR diagram: $Res^1_1,\ldots,Res^1_t$. Hence, by further  passing to a subsequence, we may assume that they all extend to homomorphisms that
factor through the reduced completion of a fixed resolution $Res^1_j$, $1 \leq j \leq t$.

By possibly passing to a further subsequence, from the weak test sequence, $\{h_n\}$ we construct a map:
$\nu:Red(Comp(Res^1_j)) \to FCl(Red(Comp(CRes)))$.`
Suppose that $FCl(Red(Comp(CRes)))$ is identical to $Red(Comp(CRes))$,
and that the map $\nu$ satisfies property (1) in lemma 1.13, i.e., suppose that $\nu$ is an isomorphism.

In case there exists such a weak test sequence $\{h_n\}$, and both maps $\nu$ and $\tau$ are isomorphisms, the map $\tau^{-1} \circ \nu$ is an
isomorphism between the reduced completions of $Res^1_j$ and $Res_1$. Hence, we can take out the resolutions $Res_1$ and $Res^1_j$ from the
two minimal complexity diagrams, and the remaining two collections of resolutions must be minimal complexity diagrams for the collection of
homomorphisms the are associated with the set $U_1$, i.e., the set of tuples from the given set $U$ that are not associated with homomorphisms that extend
to homomorphisms that factor through $Red(Comp(Res_1))$.

Therefore, we assume that there is no such weak test sequence $\{h_n\}$. In this case we look at the subset of tuples in $U$,
with associated  homomorphisms $h:F_{k+\ell} \to F_k$, that
do not extend to homomorphisms that factor through the reduced completions of $Res_2,\ldots,Res_m$. We denote this subset of tuples $U_2$.

Our goal is to show that  under the assumption that there is no weak test sequence, $\{h_n\}$, for which both maps $\nu$ and $\tau$ are isomorphisms,
it is possible to construct an MR diagram for the set $U_2$ that has strictly smaller complexity than a diagram that contains the single resolution
$Res_1$.

Such a lower complexity diagram will lead to a contradiction.
With the subset $U_2$ we can associate a minimal complexity MR diagram. Since all the tuples in $U_2$ are associated with homomorphisms that extend to
homomorphisms that factor through the reduced completion of $Res_1$, and since the minimal complexity MR diagram for the set $U_2$ can complete the resolutions
$Res_2,\ldots,Res_m$ to be an MR diagram for the set $U$, the minimal complexity MR diagram for the set $U_2$ must consist of a single resolution with
the same complexity as $Res_1$, and can not be an MR diagram with strictly smaller complexity.

Let $\{h_n:F_{k+\ell} \to F_k\}$ be a sequence of homomorphisms that are associated with tuples from the set $U_2$. By possibly passing to a
convergent subsequence, from the sequence of homomorphisms, $\{h_n\}$, we construct a resolution $CRes$. Since all the homomorphisms $\{h_n\}$ extend to
homomorphisms that factor through the reduced modeled completion, $Red(Comp(Res_1))$, by passing to a further subsequence we get a map:
$\tau: Red(Comp(Res_1)) \to FCl(Red(Comp(CRes)))$.

If $FCl(Red(Comp(CRes)))$ is not identical to $Red(Comp(CRes))$ then it has strictly smaller complexity and we continue with it. If the map
$\tau$ is not onto, we can replace the resolution $CRes$ with a resolution with strictly smaller complexity. Since complexities of resolutions
are well-ordered, after finitely many constructions of framed closures and the map $\tau$ we may assume that the map $\tau$ is onto.

If the map $\tau$ is onto, and the structure of $FCl(Red(Comp(CRes)))$ is identical to the structure $Red(Comp(CRes))$, then:
$Compl(FCl(Red(Comp(CRes)))) \, \leq \, Compl(Red(Comp(Res_1)))$. If the inequality is strict we are done. Otherwise, the map
$\tau$ must be an isomorphism by part (1) in lemma 1.13.

A subsequence of the homomorphisms, $\{h_n\}$, must factor through the reduced completion of a resolution, $Res^1_j$, from the given second minimal complexity
MR diagram. By possibly passing to a further subsequence, we get a map: $\nu:Res(Comp(Res^1_j)) \to FCl(Red(Comp(CRes)))$. If the structure of
$FCl(Red(Comp(CRes)))$ is not identical to the structure of $Red(Comp(CRes))$, or if $\tau$ is not onto,  $CRes$
can be replaced by a resolution with a strictly smaller complexity
and we are done.

Hence, by part (1) of lemma 1.13 $\nu$ is an isomorphism. But in this case there exists a sequence of homomorphisms $\{h_n\}$ that is associated with
tuples from $U_2$ for which both $\tau$ and $\nu$ are isomorphisms, which contradicts our assumptions.

Therefore, under our assumptions, from each sequence of homomorphisms that are associated with tuples from $U_2$ it is possible to extract a subsequence that
factors and is a weak test sequence w.r.t.\ a resolution $CRes$, that has complexity that is strictly bounded by the complexity of $Res_1$.

By compactness there is a finite collection of such resolutions, such that all the homomorphisms that are associated with points in $U_2$ extend to
homomorphisms that factor through their reduced completions. So we found an MR diagram for the tuples in $U_2$ that has strictly smaller
complexity than the MR diagram that contains the single resolution $Res_1$. As we already argued this contradicts the assumption the the first
given MR diagram is a minimal complexity MR diagram for the original set of tuples $U$.

Finally, this implies that the second minimal complexity MR diagram contains a resolution $Res^1_i$, for which $Red(Comp(Res^1_i))$ is
isomorphic to $Red(Comp(Res_1))$, and by continuing iteratively, the two given minimal complexity MR diagrams satisfy the conclusion of theorem
1.8.

\line{\hss$\qed$}

Theorem 1.8 proves the uniqueness of a minimal complexity MR diagram in case all the resolution in the minimal complexity diagram have reduced
completions that contain only QH vertex groups. To generalize the uniqueness of theorem 1.8 to arbitrary sets of tuples with single ended MR diagrams we need
to associate $trimmed$ reduced completion with the reduced completion of a modeled resolution.

\proclaim{Definition 1.14} Let $Res$ be a   single ended   modeled resolution, and let $Red(Comp(Res))$ be its reduced completion. In the reduced
completion. there is a canonical image of the limit group $L$, that is associated with the homomorphisms $h:F_{k+\ell} \to F_k$ from which the
completion was constructed. i.e., there is an embedding: $\rho: L \to Red(Comp(Res))$.

We set $TrimRed(Comp(Res))$, to be the subgroup (subtower) of the reduced completion that contains:
\roster
\item"{(1)}"  all the QH vertex groups, and all the edge groups along  the reduced completion $Red(Comp(Res))$.

\item"{(2)}" the image of the limit group $L$, $\rho(L)$, in $Red(Comp(Res))$.
\endroster
and where each abelian vertex group along the reduced completion is replaced with a minimal possible subgroup such that properties (1) and (2) still hold.
\endproclaim

Trimmed reduced completions of modeled resolutions enable us to generalize theorem 1.8 to sets of tuple $U$ with general minimal
complexity MR diagrams. To define MR diagrams that contains general single ended resolutions, i.e., resolutions with both QH and abelian vertex groups,
we need to slightly weaken the properties of an MR diagram that are listed in theorem 1.7.

\vglue 1.5pc
\proclaim{Definition 1.15 (cf. theorem 1.7)} Let $\ell$ be a positive integer and $U$ be a set of $\ell$-tuples. There exist finitely many (countable)
modeled resolutions such that:
\roster
\item"{(1)}" each homomorphism $h:F_{k+\ell} \to F_k$ that is associated with an $\ell$-tuple from $U$ factors through the trimmed reduced completion of
at least one of the (finitely many)
resolutions  of the (countable) modeled resolutions.

Recall (cf theorem 1.7) that
a homomorphism, $h:F_{k+\ell} \to F_k$, factors through the trimmed reduced completion of a single ended resolution $Res$, if it extends to a homomorphism of the
reduced completion, $\hat h:TrimRed(Comp(Res)) \to F_k$, and no s,c,c,\ in any of the QH vertex groups along the trimmed reduced completion is mapped by $\hat h$
to the identity in $F_k$.

\item"{(2)}" for each of the finitely many (countable) modeled resolutions there exists a sequence of homomorphisms that are associated with tuples in $U$ that form a weak
test sequence for (the completion of) the (countable) modeled resolution.
\endroster
\endproclaim

The slightly weakened properties of an MR diagram in definition 1.15 enable us to prove the uniqueness of a minimal complexity
single ended diagram that is associated with a set of tuples for general minimal complexity single ended diagrams that are associated with a set of tuples.

\vglue 1.5pc
\proclaim{Theorem 1.16 (cf. theorem 1.8)} Let $\ell$ be a positive integer and $U$ be a set of $\ell$-tuples in $F_k$. Suppose that $U$ has two minimal complexity
single ended finite collections of
(countable) modeled resolutions that satisfy the properties of definition 1.15.

Let $Res_1,\ldots,Res_m$ be the modeled resolutions in the first collection (MR diagram), and $Res^1_1,\ldots,Res^1_t$ be the modeled resolutions in the second collection.
Recall that with each resolution there is an associated embedding of a limit group $L_i$ into the reduced completion of the resolution $Res_i$, and a
similar associated embedding of a limit group $L^1_i$ into the
reduced completion of $Res^1_i$.

Then $m=t$, and up to a change of order,for each $i$, $1 \leq i \leq m$, the structures of the trimmed reduced
completions, $TrimRed(Comp(Res_i))$ and $TrimRed(Comp(Res^1_i))$, are similar (same number of levels
and same number and isomorphic  QH and abelian vertex groups  in each level), and there exist isomorphisms:
$$\eta_i:TrimRed(Comp(Res_i)) \ \to \ TrimRed(Comp(Res^1_i)) \ \ and \ \
\eta^1_i:TrimRed(Comp(Res^1_i)) \to TrimRed(Comp(Res_i))$$
that map the QH and abelian vertex groups in one completion isomorphically onto  corresponding QH and abelian vertex groups in the target completion,
and the image of the limit group $L_i$
that embeds in $TrimRed(Comp(Res_i))$ isomorphically onto the image of $L^1_i$ in $TrimRed(Comp(Res^1_i))$ and vice versa.
\endproclaim

\nfp The proof is essentially similar to the proof of theorem 1.8 with modifications to include abelian vertex groups along the resolutions.
Given a resolution from the diagram we start by trying to push down abelian and QH vertex groups using the iterative procedure that we presented. After
finitely many steps this procedure terminates.

After the procedure terminates the edge group that is connected to each abelian vertex group that is not in the bottom two levels
of the reduced modeled resolution is hyperbolic in the abelian decomposition that is associated with the next level. Each QH vertex group that is not in the bottom two levels
has a boundary subgroup that is not elliptic in the abelian decomposition that is associated with the next level.

As in the proof of theorem 1.8, and for presentation purposes, suppose first that the modeled structures of the modeled resolutions in the two minimal
complexity MR diagrams are trivial, i.e.,
that the modeled subtowers above any QH subgroup in these modeled resolutions contain only the base QH vertex groups with nothing above it. In this case
the reduced modeled resolutions are identical to the modeled resolutions.

Let $Res_1$ be a resolution of maximal complexity from the resolutions, $Res_1,\ldots,Res_m$. By the properties of the MR diagrams that are associated
with the set $U$, there exists a
sequence of homomorphisms $\{h_n\}$ that are associated with elements of the set $U$ and extend to a weak test sequence of $Res_1$.

By the properties of MR diagrams,  a subsequence of the homomorphisms, $\{h_n\}$, factor through the trimmed reduced completion of a resolution $Res^1_i$ from the
second given minimal complexity MR diagram. By the construction of
formal solutions and formal limit groups in [Se2],
from that subsequence of the sequence $\{h_n\}$ and its extensions to $Comp(Res_1)$ and to the trimmed (reduced) completion $TrimRed(Comp(Res^1_i))$
it is possible
to pass to a further subsequence that converges to a map:  $\nu:TrimRed(Comp(Res^1_i)) \to FCl(Res_1)$, where $FCl(Res_1)$ is a $framed$  closure of $Comp(Res_1)$,
and $\nu$ maps the image of the limit group $L^1_i$ in $TrimRed(Comp(Res^1_i))$ onto the image of the limit group $L$ in $FCl(Res_1)$. These enable us to
to generalize lemmas 1.9 and 1.10 to general minimal complexity single ended MR diagrams
with trivial modeled structures.

Before we examine the properties of the  map $\nu$ note that in the presence of abelian vertex groups along the completion
$Comp(Res_1)$, it may be that the structure of the trimmed completion, $TrimRed(Comp(Res_1))$, is not the same as that of $Comp(Res_1)$. i.e., it may be that
some abelian vertex groups along the completion, $Comp(Res_1)$, are replaced with abelian vertex groups of strictly smaller rank in $TrimRed(Comp(Res_1))$. If
this happens it is possible to replace the resolution $Res_1$ with a resolution with strictly smaller complexity, by replacing abelian vertex groups along
the resolution, $Res_1$, with abelian vertex groups of strictly smaller rank, and the set of homomorphisms that factor through the completion of the obtained resolution
is identical to the set of homomorphisms that factor through the completion of the original resolution, $Comp(Res_1)$. Furthermore their trimmed resolutions are isomorphic.
This clearly contradicts the assumption that the first given MR diagram is of minimal complexity.
Hence, in the sequel we may assume that the structures of $Comp(Res_1)$ and $TrimRed(Comp(Res_1))$ are similar.

\proclaim{Lemma 1.17} With the notation above if the structure of $FCl(Res_1)$ is not similar to the structure of $Comp(Res_1)$ then either:
\roster

\item"{(1)}" two non-conjugate abelian or   QH vertex groups in $Comp(Res_1)$ are mapped into conjugate ones
in $FCl(Res_1)$.

\item"{(2)}" Some QH vertex groups in the framed closure $FCl(Res_1)$ properly contains a  QH vertex groups in  $Comp(Res_1)$.
\endroster
\endproclaim

\nfp Similar to lemma 1.9 and follows from the construction of a closure in [Se2].

\line{\hss$\qed$}

As in lemma 1.9, in both cases in lemma 1.17 the complexity of $FCl(Res_1)$ is strictly smaller than the complexity of $Comp(Res_1)$.
Now, suppose that the structures of
$FCl(Res_1)$ and $Comp(Res_1)$ are similar.

\proclaim{Lemma 1.18} If $FCl(Res_1)$ and $Comp(Res_1)$ have similar structure then either:
\roster
\item"{(1)}" The restricted map $\nu:TrimRed(Comp(Res^1_i)) \to FCl(Res_1)$  is an isomorphism from $TrimRed(Comp(Res^1_i))$ onto
$TrimRed(Comp(Res_1))$. In particular, the two completions $TrimRed(Comp(Res_1))$ and $TrimRed(Comp(Res^1_i)|)$ have the same
structure and the homomorphisms $\{h\}$ that are associated with tuples from $U$ and extend to homomorphisms that factor through $TrimRed(Comp(Res_1))$ extend to
homomorphisms that factor through $TrimRed(Comp(Res^1_i))$ and vice versa.

\item"{(2)}"  The  map $\nu$ does not satisfy part (1). In that case either there exists a QH vertex group $Q$ in $FCl(Res_1)$ that the image
$\nu(TrimRed(Comp(Res^1_i)))$ intersect
trivially all its conjugates, i.e., it intersects all its conjugates only in conjugates of some boundary subgroups, or there is an abelian vertex group $A$ in
$FCl(Res_1)$ that $\nu(TrimRed(Comp(Res^1_i)))$ intersects all the conjugates of $A$ in subgroups that are contained in  conjugates of a subgroup $\hat A < A$, and
$\hat A$ is a subgroup of infinite index in $A$.
\endroster
\endproclaim

\nfp
The argument that we use to prove lemma 1.18 is a modification of the proof of lemma 1.10.

If the completion of $Res_1$ contains
no abelian vertex groups the lemma follows from lemma 1.10. Hence, we can assume that $Res_1$ contains abelian vertex groups.
Let $Q_1$ be a QH vertex group in the terminal  abelian decomposition  of $Comp(Res^1_i)$. All the boundary components of $Q_1$ can be conjugated
into the terminal level of $Comp(Res^1_i)$.

We treat $Q_1$ as we treated QH vertex groups in the terminal abelian decomposition of $Comp(Res^1_i)$ in the proof of lemma 1.10 (in case there are
no abelian vertex groups in $Comp(Res_1)$). If $Q_1$ is not mapped isomorphically onto a  surviving QH vertex group in $FCl(Res_1)$, by a composition of $\nu$ with appropriate
retractions of $FCl(Res_1)$ (see definition 1.11), then all the union of QH vertex groups in $FCl(Res_1)$ that $\nu(Q_1)$  intersects non-trivially a conjugate of them have
strictly smaller complexity than $Q_1$. If $Q_1$ is mapped isomorphically onto a conjugate of a surviving surface in $FCl(Res_1)$,  then we modify the map $\nu$
and $FCl(Res_1)$ so that $Q_1$ is mapped isomorphically onto a QH vertex group
in the terminal abelian decomposition of $FCl(Res_1)$.

We continue by analyzing the images of all the QH vertex groups in the terminal abelian decomposition of $Comp(Res^1_i)$. Let $M_1$ be the subgroup of
$Comp(Res^1_i)$ that contains the terminal vertex group in $Comp(Res^1_i)$ and all the QH vertex groups in the terminal abelian decomposition in $Comp(Res^1_i)$.
 Suppose that at least one of the QH vertex groups in $M_1$
 is not mapped isomorphically onto a conjugate of a surviving surface in $FCl(Res_1)$ (by a composition of $\nu$ with appropriate retractions). Then
the total complexity of the QH vertex groups in $Comp(Res_1)$ that $\nu(M_1)$ intersects a conjugate of them non-trivially is strictly
smaller than the complexity of the total complexity of the QH vertex groups in the terminal abelian decomposition in $Comp(Res^1_i)$.

By climbing along the levels of $Comp(Res^1_i)$ as we did in the proof of lemma 1.10, it follows that in this case $\nu(TrimRed(Comp(Res^1_i)))$ intersects trivially
the conjugates of at least one QH vertex group in $FCl(Res_1)$. Hence, in this case part (2) of the lemma follows.

Therefore, we may assume that all the QH vertex groups in the terminal abelian decomposition of $Comp(Res^1_i)$ are mapped isomorphically onto
conjugates of
surviving QH vertex groups in $FCl(Res_1)$.  In that case we continue with the abelian vertex groups in the terminal abelian decomposition in $TrimRed(Comp(Res^1_i))$.

The edge group that is connected to such an abelian vertex group $A_1$ is contained in the terminal group is $TrimRed(Comp(Res^1_i))$ that is mapped by $\nu$ to the terminal vertex group
in $FCl(Res_1)$. Hence, $\nu(A_1)$ is either contained in the terminal limit group of $FCl(Res_1)$ or it is contained in a conjugate of an abelian vertex group in the
terminal abelian decomposition of $FCl(Res_1)$.

As in the case of a QH vertex group, we say that an abelian vertex group $A$ in the terminal abelian decomposition of $FCl(Res_1)$ is $surviving$ if there exists
an abelian vertex group $A_1$ in $TrimRed(Comp(Res^1_i))$, such that:
\roster
\item"{(1)}"  $\nu$ maps $A_1$ isomorphically onto a finite index subgroup of a conjugate of $A$.

\item"{(2)}" $\nu$ maps the other abelian vertex groups in the terminal abelian decomposition of $TrimRed(Comp(Res^1_i))$ to subgroups that intersect conjugates
of $A$ trivially.
\endroster

As for  QH vertex group in the terminal abelian decomposition of $Comp(Res^1_i)$, suppose that
 there exists an abelian vertex group $A_1$ in that terminal abelian decomposition of $TrimRed(Comp(Res^1_i))$
that is not mapped isomorphically onto a finite index subgroup of
a conjugate of a surviving abelian vertex group in
the terminal abelian decomposition of $FCl(Res_1)$.

In that case we continue to climb along the levels of $TrimRed(Comp(Res^1_i))$ and $FCl(Res_1)$, as we did in analyzing the map $\nu$ in completions that
don't have abelian vertex groups (lemma 1.10), and either:
\roster
\item"{(1)}"  there exists a QH vertex group in one of the levels
of $FCl(Res_1)$ that all its conjugates intersects trivially the image, $\nu(TrimRed(Comp(Res^1_i)))$.

\item"{(2)}" there exists an abelian vertex group $A$ along the levels of $FCl(Res_1)$ that $\nu(TrimRed(Comp(Res_1)))$ intersects all the conjugates of $A$ in subgroups
that are contained in conjugates of a subgroup of infinite index $\hat A<A$.
\endroster

In both  cases part (2) of the lemma holds. Hence, we may assume that after modifying the map $\nu$ and $FCl(Res_1)$ as we did in the proof of lemma 1.10,
 all the QH vertex groups in the terminal abelian decomposition of $TrimRed(Comp(Res_1))$
are mapped by $\nu$ onto conjugates of surviving QH vertex groups in the terminal abelian decomposition of $FCl(Res_1)$,
and all the abelian vertex groups in the terminal abelian decomposition of $TrimRed(Comp(Res_1))$
are mapped isomorphically onto finite index subgroups of conjugates of surviving abelian vertex groups in the terminal abelian decomposition of $FCl(Res_1)$.

\smallskip
We continue iteratively by climbing along the levels of $TrimRed(Comp(Res^1_i))$,
and possibly further modifying the structure of
$FCl(Res_1)$ and the map $\nu$ precisely as we did in the proof of lemma 1.10. In each level in $TrimRed(Comp(Res^1_i))$ either all the QH  vertex groups are mapped
isomorphically onto surviving QH vertex groups  in $FCl(Res_1)$, and all the abelian vertex groups are mapped isomorphically onto finite index subgroups
of surviving abelian vertex groups in $FCl(Res_1)$, or part (2) of  lemma 1.18 holds.

The analysis is complete when we analyze the highest level in $TrimRed(Comp(Res^1_i))$. The two given MR diagrams were assumed to be of minimal
complexity, and the resolution $Res_1$ is assumed to be of maximal complexity among the resolutions in the first diagram. Hence,  if part (2) of lemma 1.18
does not hold, then $\nu$ maps isomorphically every QH vertex in $TrimRed(Comp(Res^1_i))$ onto a conjugate of a QH vertex group in the same level $FCl(Res_1)$,
and every abelian  vertex group in $TrimRed(Comp(Res^1_i))$ isomorphically onto a conjugate of a finite index subgroup of an abelian vertex group
in the same level in $FCl(Res_1)$.

Suppose that part (2) in lemma 1.18 does not hold. In that case each abelian vertex group $A_1$  in $TrimRed(Comp(Res^1_i))$ is mapped isomorphically onto a finite
index subgroup of a conjugate of an abelian vertex $A$ in the same level of $FCl(Res_1)$. By the construction of the trimmed reduced
modeled completion,  $TrimRed(Comp(Res^i_1))$, which is a subgroup of $Comp(Res^1_i)$, contains the image of  $L^1_i$ in $Comp(Res^1_i)$
 as well as all the edge groups that are connected to QH vertex groups and to abelian vertex groups along the levels of $Comp(Res^1_i)$.

Since part (2) of lemma 1.18 does not hold, $\nu$ is an isomorphism that maps QH vertex groups in $TrimRed(Comp(Res^1_i))$ isomorphically onto conjugates of QH
vertex groups in $Comp(Res_1)$, and abelian vertex groups in $TrimRed(Comp(Res^1_i))$ isomorphically onto finite index subgroups of conjugates of
abelian vertex groups in $Comp(Res_1)$. Hence, $\nu(TrimRed(Comp(Res^1_i))$ contains the image of $L_1$ in $Comp(Res_1)$ as well as conjugates of all the
edge groups that are connected to abelian and QH vertex groups along the levels of $Comp(Res_1)$. Furthermore, $\nu(TrimRed(Comp(Res^1_i)))$ is obtained
from $FCl(Res_1)$ by replacing each of the abelian vertex groups along $FCl(Res_1)$ with a subgroup of finite index, which is conjugate of the image
of a corresponding abelian vertex group in $TrimRed(Comp(Res^1_i))$.

Let $A$ be an abelian vertex group in $FCL(Res_1)$, and let $TrimRed(A)$ be the finite index subgroup of $A$ that replaces $A$ in the trimmed reduced
completion of $Comp(Res_1)$, $TrimRed(Comp(Res_1))$. Let $A_1$ be a vertex group in $TrimRed(Comp(Res^1_i))$ that is mapped by $\nu$ onto a conjugate
of a finite index subgroup of $A$.

If $TrimRed(A)$ is not conjugate to $\nu(A_1)$ in $FCl(Res_1)$, then $TrimRed(A)$ intersects a conjugate of $\nu(A_1)$ in a proper finite index subgroup of
either $TrimRed(A)$ or of a conjugate of $\nu(A_1)$. If it is a proper subgroup of $TrimRed(A)$, then it is possible to replace $TrimRed(A)$ be the
proper intersection and obtain a proper subgroup of the trimmed reduced completions that has the same structure of $TrimRed(Comp(Res_1))$, and still
contains the image of $L_1$ in $Comp(Res_1)$ and conjugates of all the edge groups along the levels of $Comp(Res_1)$, a contradiction to the minimality of
the trimmed reduced completion, $TrimRed(Comp(Res_1))$.

If the intersection of $TrimResd(A)$ and a conjugate of $\nu(A_1)$ is a proper subgroup of $\nu(A_1)$, then it is possible to replace $TrimRed(Comp(Res^1_i))$
with a proper subgroup that contains the image of $L^1_i$ in $Comp(Res^1_i)$ and conjugates of all the edge groups along the levels of $Comp(Res^1_i)$,
a contradiction to the minimality of $TrimRed(Comp(Res^1_i))$.

Therefore, $\nu$ must map $TrimRed(Comp(Res^1_i))$ isomorphically onto $TrimRed(Comp(Res_1))$ in $FCl(Res_1)$, and case (1) of lemma 1.18 follows.

\line{\hss$\qed$}

Lemmas 1.17 and 1.18, like lemmas 1.9 and 1.10, were stated and proved in the case in which the modeled completions of the resolutions in the minimal complexity
diagrams are trivial, i.e., in the case in which the reduced model completions of the resolutions is identical to the modeled completions. Like in the case in which
there are no abelian vertex groups, in which lemmas 1.9 and 1.10 were generalized in lemmas 1.12 and 1.13 to general single ended resolutions that contain
no abelian vertex groups, and by the same argument,
lemmas 1.17 and 1.18 remain valid for general single ended resolutions, i.e. to single ended resolutions with possibly non-trivial modeled structures, which means
that they may have reduced completions that  defer from
the completions.

Given lemmas 1.17 and 1.18 for general single ended completions, the proof of theorem 1.16 using these lemmas is identical to the proof of theorem 1.8 from
lemmas 1.12 and 1.13.

\line{\hss$\qed$}

Let $\ell$ be a positive integer and $U$ be a set of $\ell$-tuples in $F_k$. Suppose that it is possible to associate
with $U$ a minimal complexity
single ended MR diagram.

By theorem 1.16 it is possible to associate with $U$ a canonical minimal complexity diagram, and this minimal complexity diagram is unique up to
isomorphisms of the trimmed completions of the resolutions in the diagram.

The tools and the arguments that are used in the proof of theorem 1.16 generalize to torsion-free hyperbolic groups using the constructions and the results in
[Se6]. Hence, a canonical single ended diagram can be associated with sets of tuples in a torsion-free hyperbolic group, in the same way it is
associated with sets of tuples over a free group.

\vglue 1.5pc
\proclaim{Theorem 1.19 (cf. theorem 1.16)} Let $\Gamma$ be a torsion-free hyperbolic group, let $\ell$ be a positive integer, and let
$U$ be a set of $\ell$-tuples in $\Gamma$. Suppose that it is possible to associate with $U$ a single ended MR diagram.

Then it is possible to associate with $U$ a canonical minimal complexity diagram, and this minimal complexity diagram is unique up to
isomorphisms of the trimmed completions of the resolutions in the diagram.
\endproclaim

Note that theorem 1.19 is stated for torsion-free hyperbolic groups. To generalize theorem 1.19 to hyperbolic groups with torsion, part of the needed
tools and techniques appear in [Re-We] and [He], and the remaining ones will probably be available in a current ongoing work of Gili Berk.

\vglue 1.5pc
\centerline{\bf{\S2.  A pseudo topology on ordered  tuples, their rank,  }}
\centerline{\bf{  and the canonical envelope of a definable set}}
\medskip
\medskip

Let $F_k$ be a non-abelian free group, and let $\ell$ be a positive integer. In the previous section we have associated
a canonical (minimal complexity) MR diagram with any given set $S \subset F_k^{\ell}$ in case it is possible to associate a single ended diagram with $S$ (theorem 1.16).

The canonical diagram certainly associates a complexity with the set $S$, and it also associates a $closure$ of such a set $S$, although the set of closures does not define a topology,
but rather a $pseudo$ topology.

\vglue 1.5pc
\proclaim{Definition 2.1} Let $F_k$ be free group, and $\ell$ a positive integer. Let $S \subset {F_k}^{\ell}$ be a set of (ordered) $\ell$-tuples with a single ended canonical
(minimal complexity) diagram $\Delta$ (by theorem 1.16). By definition, every tuple from $S$ corresponds to a homomorphism, $h:F_{k+\ell} \to F_k$,
and every such homomorphism extends to a homomorphism that factors through at least  one of the trimmed reduced completions of the resolutions in $\Delta$.

We define the pseudo-closure of $S$, $pcl(S)$, to be the set of all the $\ell$-tuples from $F_k$  that correspond to homomorphisms $f:F_{k+\ell} \to F_k$,
that extend to
homomorphisms, $\hat f:TrimRed(Comp(Res))$, for some of the trimmed reduced completions in the canonical MR diagram $\Delta$ of $S$, so that
$\hat h$ maps no s.c.c.\ on any of the QH vertex groups in $TrimRed(Comp(Res))$ to the identity in $F_k$.
\endproclaim

Note that by definition $\Delta$, which is the canonical diagram of $S$ is also the  canonical diagram of $pcl(S)$, so $pcl(pcl(S))=S$. Unfortunately,
the intersection and the union of two closed sets
may not be closed.
However, we have the following.

\vglue 1.5pc
\proclaim{Theorem 2.2} Let $S_1,S_2$ satisfy the assumptions on sets of $\ell$-tuples in definition 2.1, and let $\Delta_1,Delta_2$ be their canonical diagrams. Suppose that
the complexity of $\Delta_1$ is bounded by the complexity of $\Delta_2$. Then:
\roster
\item"{(1)}" $pcl(S_1) \, \cap \, pcl(S_2)$ has a canonical diagram $\Delta$ that is either $\Delta_1$ or  has strictly smaller complexity than $\Delta_1$.

\item"{(2)}" $pcl(S_1) \, \cup \, pcl(S_2)$ has a canonical diagram $\Delta'$, that is either of the union of the resolutions in $\Delta_1$ and $\Delta_2$, or  its
complexity is strictly bounded by the complexity of $\Delta_1 \cup \Delta_2$. Furthermore, all the maximal complexity resolutions in $\Delta_1 \cup \Delta_2$ are also the maximal
complexity resolutions in $\Delta'$.

\item"{(3)}" If $S_1 \subset S_2$ then  necessarily the complexity of $\Delta_1$ is bounded above by the complexity of $\Delta_2$, and if the complexities are
equal  then
$Delta_1=\Delta_2$ (but it may be that $pcl(S_1)$ is not contained
in $pcl(S_2)$).

\item"{(4)}" every variety  over a free group with a single ended MR diagram is pseudo-closed.
\endroster
\endproclaim

\nfp Let $T=pcl(S_1) \cap pcl(S_2)$. The canonical diagrams $\Delta_1$ and $\Delta_2$ are MR diagrams for the intersection $T$, and both are assumed to be
single ended. Hence, $T$ has a single ended MR diagram, and by theorem 1.16 $T$ has a canonical MR diagram that we denote $\Delta$. Since $\Delta$ is
of minimal complexity, and $\Delta_1$ and $\Delta_2$ are MR diagrams for $T$, the complexity of $\Delta$ is bounded above by the complexities of
$\Delta_1$ and $\Delta_2$. Hence, since the complexity of $\Delta_1$ is bounded above by the complexity of $\Delta_2$, by theorem 1.16, either $\Delta$ is
equal to $\Delta_1$ (up to isomorphism) or the complexities of both $\Delta_1$ and $\Delta_2$ are strictly bigger than the complexity of $\Delta$.

To prove (2) let $R=pcl(S_1) \cup pcl(S_2)$. The union of the resolutions in both $\Delta_1$ and $\Delta_2$ is an MR diagram for $R$, so $R$ has a single ended MR diagram.
Hence, $R$ has a canonical minimal complexity diagram by theorem 1.16. Let $\Delta'$ be the minimal complexity (canonical) MR diagram for the set $R$. Since the union of the
resolutions in $\Delta_1$ and $\Delta_2$ is an MR diagram for $R$, the complexity of $\Delta'$ is bounded above by the complexity of $\Delta_1 \cup \Delta_2$. By theorem
1.16 if the complexities of $Delta'$ and $\Delta_1 \cup \Delta_2$ are equal, then the trimmed reduced completions of the resolutions in the two diagrams are isomorphic.

Let $Res$ be a maximal complexity resolution in $\Delta_1 \cup \Delta_2$. In $Comp(Res)$ there is an embedded limit group $L$. Every $\ell$-tuple
that corresponds to a homomorphism of the
limit group $L$ into $F_k$ that extends to a homomorphism of $TrimRed(Comp(Res))$ is in the set $R$.

We look at all the weak test sequences of the trimmed reduced completion, $TrimRed(Comp(Res))$.
Since the restriction of every homomorphism of $TrimRed(Comp(Res))$
to the limit group $L$ is in $R$, the restrictions of a test sequence  of the trimmed reduced completion of $Res$ must extend to homomorphisms that factor through the
trimmed completion of a resolution in $\Delta'$. Hence, using the techniques of [Se2] that were applied in the previous section, given a test sequence of
$TrimRed(Comp(Res))$ we can pass to a subsequence, from which it is possible to get a map: $\nu:TrimRed(Comp(Res')) \to FCl(TrimRed(Comp(Res)))$, where
$Res'$ is a resolution from the MR diagram $\Delta'$.

If for every test sequence of homomorphisms of $TrimRed(Comp(Res))$, the structure of $FCl(TrimRed(Comp(Res)))$ is not the same as the structure
of $TrimRed(Comp(Res))$, or their structures are the same, but the map $\nu$ satisfies part (2) in lemma 1.18, then the maximal complexity resolution
$Res$ in either $\Delta_1$ or $\Delta_2$ can be replaced by finitely many resolutions of strictly smaller complexity. This contradicts the minimal
complexity of either $\Delta_1$ or $\Delta_2$ (for $S_1$ or $S_2$ in correspondence).

Hence, there exists a test sequence of $TrimRed(Comp(Res))$ for which   $FCl(TrimRed(Comp(Res)))$ has the same structure as $TrimRed(Comp(Res))$,
and the map $\nu$ does not satisfy condition (2) in lemma 1.18.
In that case
the complexity of $TrimRed(Comp(Res'))$ is bounded below by the complexity of $TrimRed(Comp(Res))$.

Since $\Delta'$ is a minimal complexity MR diagram for $R$, its
complexity is bounded by the complexity of $\Delta_1 \cup \Delta_2$. $Res$ is the maximal complexity resolution in $\Delta_1 \cup \Delta_2$ so its complexity
bounds from above the  complexities of the  resolutions in $\Delta'$, $Res'$ is a resolution in $\Delta'$, so the complexity of $TrimRed(Comp(Res))$ bounds from
above the complexity of $TrimRed(Comp(Res'))$. Hence, the complexities of $TrimRed(Comp(Res'))$ and $TrimRed(Comp(Res_1))$ are equal. Therefore, by part (1) of lemma 1.18,
$\nu$ restricts to an isomorphism between $TrimRed(Comp(Res'))$ and $TrimRed(Comp(Res))$.

Therefore, for each maximal complexity resolution $Res$ in $\Delta_1 \cup \Delta_2$ there exists a maximal complexity resolution $Res'$ in $\Delta'$ such that
$Res$ and $Res'$ have  isomorphic
trimmed reduced completions. So the set of maximal complexity resolutions in $\Delta'$ and $\Delta_1 \cup \Delta_2$ is similar up to isomorphism (this need not be the case
for lower complexity resolutions in the two diagrams).

To prove (3) note that if $S_1 \subset S_2$ then every MR diagram for $S_2$ is an MR diagram for $S_1$. Hence, if $S_2$ has a single ended MR diagram so is $S_1$, and the complexity
of the canonical MR diagram  for $S_2$ bounds the complexity of the canonical (minimal complexity) MR diagram for $S_1$. If their complexities are equal, then by theorem
1.16 the two diagrams are essentially isomorphic.

To prove (4) let $S$ be a variety over $F_k$ with a single ended MR diagram. By theorem 1.16 $S$ has a canonical (minimal complexity) diagram
that we denote $\Delta$.
If $Res$ is a resolutions in $\Delta$, then by construction there is a map of a limit
group $L$ into $Comp(Res)$. Since, by construction (see the proof of theorem 1.16)  there exists a weak test sequence of homomorphisms of $Res$, the limit group $L$ satisfies the
equations of the variety that is associated with $S$. Hence, all the specializations of $L$ that extend to homomorphisms that factor through $TrimRed(Comp(Res))$ must satisfy all
the equations that $L$ satisfies, so they are all in the variety $L$. Therefore, $S=pcl(S)$.

\line{\hss$\qed$}

The Zariski topology over free and hyperbolic groups is Noetherian since these groups are equationally Noetherian.
Noetherianity holds for the pseudo topology that we introduced,
because of the properties of the complexities of single ended MR diagrams.

\vglue 1.5pc
\proclaim{Proposition 2.3} Let $C_1,\ldots$ be a family of pseudo closed sets. Let $I_j=pcl(C_1 \cap \ldots \cap C_j)$, and suppose that $I_j \neq I_{j+1}$.
Then the sequence terminates after finitely many steps.
\endproclaim

\nfp By construction the complexities of the pseudo-closed sets $I_j$ is strictly decreasing. Since the complexities of single ended diagrams are a
well-ordered set, every strictly decreasing set of complexities terminates. Hence, the sequence $I_j$ terminates after finitely many steps.

\line{\hss$\qed$}

The Noetherianity of pseudo closed sets, and their complexities enable us to define the $rank$ of a set of tuples,
which is a countable ordinal that one can associate with such a set, in a similar way that one assign Shelah, Lascar and Morely ranks with definable types
and sets in superstable theories in model theory.

\vglue 1.5pc
\proclaim{Definition 2.4} Let $\ell>0$ and let $S$ be a set of $\ell$-tuples from a free group $F_k$. Suppose that $S$ has a single ended MR diagram, hence, by
theorem 1.16 $S$ has a canonical MR diagram $\Delta(S)$. Note that this implies that every subset $\hat S \subset S$ has a canonical (single ended) diagram as well,
and its complexity is bounded by the complexity of $\Delta(S)$.

With $S$ we associate a countable ordinal, $\alpha(S)$,  that we view as its $rank$. We define the rank inductively (cf. section 8.6 in [Te-Zi] and
section 17 in [Po2]).
\roster
\item"{(1)}" $\alpha(S) \geq 0$.

\item"{(2)}" $\alpha(S) \geq \beta+1$ if there exists a subset $\hat S \subset S$, with a canonical diagram $\Delta(\hat S)$, such that the
complexity of $\Delta(\hat S)$ is strictly smaller than the complexity of $\Delta(S)$, and $\alpha(\hat S)=\beta$.

\item"{(3)}" for a limit ordinal $\lambda$, $\alpha(S) \geq \lambda$ if for every $\beta < \lambda$, there exists a subset $\hat S \subset S$, such that the
complexity of $\Delta(\hat S)$ is strictly smaller than the complexity of $\Delta(S)$, and $\alpha(\hat S)=\beta$.
\endroster
\endproclaim

By the results of [Se6] all the conclusions of theorem 2.2  and proposition 2.3 are valid for torsion-free hyperbolic groups, as well as the rank of a set that is
defined in definition 2.4.
In the second section in [Se8] we introduced
an $envelope$ of a definable set, and later used it to prove that certain sets are not definable, to classify definable equivalence relations over free and hyperbolic groups, and explicitly
to prove the stability of free and hyperbolic groups in [Se7].

The construction of the envelope in [Se8] proves its existence, but it doesn't show that the envelope is canonical. The construction in [Se8] starts with a bigger set that contains the definable set,
e.g., it's Zariski closure, and iteratively reduces it using the sieve procedure from [Se5], until an envelope is constructed.
The canonical MR diagram that is associated with a set, and the pseudo-closure of it,
enable one to associate a canonical (minimal complexity) envelope with a minimal rank definable set.

\vglue 1.5pc
\proclaim{Theorem 2.5} With the assumptions and the notation of theorem 2.1 let $S$ be a definable set with a single ended MR diagram
(over a free or a torsion-free hyperbolic group).
Then $S$ has a canonical minimal complexity envelope that is unique up to
an isomorphism of the trimmed reduced completions of the resolutions in the canonical envelope.
\endproclaim

\nfp Let $S$ be a definable set with a single ended MR diagram. By theorem 1.16 $S$ has a canonical (minimal complexity) MR diagram $\Delta$. Let $Res$
be a resolution in $\Delta$.

Given the resolution $Res$ and its trimmed reduced completion, $TrimRed(Comp(Res))$, we apply the construction of the envelope of a definable set
that appears in section1 in [Se8], and associate with it a finite set of $closures$, $Cl(TrimRed(Comp(Res)))$, that determine which cosets of specializations
of $TrimRed(Comp(Res))$ contain test sequences that are in the definable set $S$ (note that test sequences were defined in section 1 in [Se2] and the construction
of the closures and the cosets that contain test sequences in the definable set in section 1 of [Se8]).

We continue with the set of closures that are associated with  $TrimRed(Comp(Res))$ for all the resolutions $Res$ in $\Delta$.
We further continue with those elements in $S$ that do not extend to homomorphisms in any of the cosets that are associated with  $TrimRed(Comp(Res))$ and contain
test sequences of specializations that restrict to elements in $S$. We denote this set $S_1$.

$S_1$ being a subset of $S$ has a single ended MR diagram, hence, by theorem 1.16, it has a canonical MR diagram $\Delta_1$. $\Delta_1$ has strictly smaller
complexity than $\Delta$, and by construction it does not contain any of the resolutions in $\Delta$. With the resolutions in $\Delta_1$ we associate a finite
set of closures as we did with the resolutions in $\Delta$.

We define $S_2$ to be the set of elements in $S_1$ that do not extend to homomorphisms in any of the cosets that are associated with the trimmed reduced
completions of resolutions in $\Delta_1$ and contain tests sequences of specializations that restrict to elements in $S_1$.

We continue iteratively. Since at each step the complexity of the constructed diagram is strictly smaller than the complexity of the previous
diagram, and the set of complexities is well-ordered, the construction terminates after finitely many steps. Finally, the envelope of the definable set $S$
is the union of the resolutions in all the constructed diagrams together with the sets of closures that were associated with each resolution.

\line{\hss$\qed$}

\vglue 1.5pc
\centerline{\bf{\S3.  Sets of tuples in a free semigroup}}
\medskip

In the first section we associated a canonical MR diagram with a set of tuples in a free (or a torsion-free hyperbolic), and in the second section
we used the diagram to define a pseudo-topology on the set of tuples. In this section we prove analogous results for sets of tuples in a free semigroup.

Let $FS_k$ be a free semigroup of rank $k$, and let $\ell$ be a positive integer. In [Se9] we constructed a Makanin-Raborov diagram for a variety over
a free semigroup. A variety over a free semigroup, i.e., the set of solutions to a system of equations over a free semigroup,
corresponds to all the homomorphisms $\{h: SG \to FS_k\}$, where $SG$ is a f.g.\ semigroup. Since a free semigroup is embedded in the free group,
every such homomorphism $h$ extends to a pair homomorphism: $hp:(SG,G) \to (FS_k,F_k)$, where $G$ is a f.g.\ group with the same presentation as the semigroup $SG$, but the
semigroup relations of $SG$ are interpreted as group relations in $G$, and the image of the semigroup $SG$ generates $G$ as a group.
In [Se9] we have shown how to associate with the collection of pair homomorphisms $hp:((SG,G) \to (FS_k,F_k)$ a (non-canonical)  Makanin-Razborov diagram.

Let $U$ be a set of $\ell$-tuples in the free semigroup $FS_k$. As in groups, each $\ell$-tuple can be viewed as a homomorphism: $h:FS_{k+\ell} \to FS_k$. As we did in the analysis of
varieties, every such semigroup homomorphism $h$ extends to a pair homomorphism $ph:(FS_{k+\ell},F_{k+\ell})\to (FS_k,F_k)$.

In section 1 we constructed an MR diagram from a given collection of group homomorphisms that are associated with a set of tuples (theorem 1.2). With each resolution in the constructed
MR diagram there is a map from the free group $F_k$ into the reduced completion of the resolution, where the image of the map is a limit group. Using the techniques and results of [Se9],
from a given collection of pair homomorphisms, that are associated with a set of $\ell$-tuples $U \subset FS_k$, it is possible to construct an MR diagram. With each resolution in such
a diagram there is a map of pairs from $(FS_{k+\ell},F_{k+\ell})$ into the reduced completion of the resolution, where the image of the pair homomorphisms  is a
limit pair $(S,L)$. The obtained MR diagram satisfy similar properties as those that are listed in theorem 1.2.

\vglue 1.5pc
\proclaim{Theorem 3.1} Let $FS_k$ be a free semigroup,  $\ell$ be a positive integer, and $U$  a set of $\ell$-tuples in $FS_k$. There exists finitely many resolutions such that:
\roster
\item"{(1)}" all the pair homomorphisms $ph:(FS_{k+\ell},F_{k+\ell}) \to (FS_k,F_k)$ that are associated with $\ell$-tuples from $U$ factor through at least one of the finitely many
resolutions.

\item"{(2)}" for each of the finitely many resolutions, there exists a sequence of pair homomorphisms that are associated with tuples in $U$ that converge into
the limit groups along the resolutions and to their associated abelian decompositions (after finitely many  refinings).

\item"{(3)}" with each of the finitely many resolutions there is an associated pair map: $\rho_i: (S_i,L_i) \to Comp(Res_i)$, where $(S_i,L_i)$ is the limit of the original
sequence of pair homomorphisms $\{h_n\}$ that were used to construct the resolution.
\endroster
\endproclaim

Starting with the construction of some MR diagram for a given set of $\ell$-tuples, and assuming the MR diagram is single ended, we defined a complexity of MR diagrams and proved
that minimal complexity single ended MR diagram that is associated with a given set of $\ell$-tuples from a free group is unique up to isomorphism of trimmed reduced completions of
the resolutions in the diagram (theorem 1.16).

Applying these arguments for a set of pair homomorphisms that are associated with a given set of $\ell$-tuples in a free semigroup, we get a similar result for sets
of $\ell$-tuples in a free semigroup.

\vglue 1.5pc
\proclaim{Theorem 3.2 (cf. theorem 1.16)} Let $\ell$ be a positive integer and $U$ be a set of $\ell$-tuples in $FS_k$. Suppose that with $U$
it is possible to associate a single ended MR diagram that satisfies the properties that are listed in theorem 3.1.
Then the trimmed reduced completions of
the resolutions in a minimal complexity MR diagram that is associated with $U$ are unique up to isomorphism.

Suppose that $U$ has two minimal complexity
single ended finite collections of
(countable) modeled resolutions that form an MR diagram for (the pair homomorphisms that are associated with) the set of $\ell$-tuples $U$.

Let $Res_1,\ldots,Res_m$ be the modeled resolutions in the first collection (MR diagram), and $Res^1_1,\ldots,Res^1_t$ be the modeled resolutions in the second collection.
Recall that with each resolution there is an associated embedding of a limit pair $(S_i,L_i)$ into the reduced completion of the resolution $Res_i$, and a
similar associated embedding of a limit pair $(S^1_i,L^1_i)$ into the
reduced completion of $Res^1_i$.

Then $m=t$, and up to a change of order,for each $i$, $1 \leq i \leq m$, the structures of the trimmed reduced
completions, $TrimRed(Comp(Res_i))$ and $TrimRed(Comp(Res^1_i))$, are similar,
and there exist isomorphisms:
$$\eta_i:TrimRed(Comp(Res_i)) \ \to \ TrimRed(Comp(Res^1_i)) \ \ and \ \
\eta^1_i:TrimRed(Comp(Res^1_i)) \to TrimRed(Comp(Res_i))$$
that map the QH and abelian vertex groups in one completion isomorphically onto  corresponding QH and abelian vertex groups in the target completion,
and the image of the limit pair $(S_i,L_i)$
that embeds in $TrimRed(Comp(Res_i))$ isomorphically onto the image of $(S^1_i,L^1_i)$ in $TrimRed(Comp(Res^1_i))$ and vice versa.
\endproclaim

\nfp Identical to the proof of theorem 1.16.

\line{\hss$\qed$}

Given a canonical MR diagram that is associated with a set of $\ell$-tuples in a free semigroup, $FS_k$, we can define the pseudo closure of a set of tuples
if this canonical diagram is single ended. Since the construction of the canonical diagram and its properties are similar to the diagram that is associated with
a set of $\ell$-tuples in a free group $F_k$, the properties of pseudo closures of sets of tuples in a free (or hyperbolic) group that are listed in theorem 2.2
and proposition 2.3 are valid to closures of sets of $\ell$-tuples in a free semigroup as well.

\vglue 1.5pc
\centerline{\bf{\S4.  Sets of homogeneous tuples in a free associative algebra}}
\medskip

In the first section we associated a canonical MR diagram with a set of tuples in a free (or a torsion-free hyperbolic), and in the third section
we used the same construction to associate a canonical MR diagram with a set of tuples in a free semigroup. In this section we use the same techniques to associate
a canonical diagram with a set of homogeneous tuples in a free associative algebra.

Let $F$ be a finite field and let
$FA_k$ be the free associative algebra of degree $k$ over $F$. Let $\ell$ be a positive integer, and let $S$ be a set of ordered $\ell$-tuples of homogeneous elements from $FA_k$.
With each element from
$S$ we can naturally associate a homomorphism $\hat h:FA_{\ell} \to FA_k$.

Let $w \in FA_k$ be  a homogeneous element. Any two presentations of $w$ as a product of non-trivial homogeneous elements have a common refinement. Hence,
given  $w \in FA_k$, there exists the most refined presentation of $w$ as such a product.

We further set $U$ to be the set of  all the homogeneous elements in $FA_k$ that can not be
presented as a non-trivial product of homogeneous elements. The semigroup that is generated by all the elements in $U$ is $FS(U)$, the free semigroup
that is freely generated by the elements in $U$.

Every homogeneous element $w \in FA_k$ can be presented uniquely as a word of elements in $U$, so with any
homogeneous element $w \in F_k$ we can assign an element in $FS(U)$. Hence, with any $\ell$-tuple
of homogeneous elements in $FA_k$ we can associate a homomorphism $h:FS_{\ell} \to FS(U)$.

$FS(U)$ is not finitely but rather countably generated. Hence, we can apply the techniques that we used in section 1 to construct MR diagrams from
subsequences of homomorphisms but
some modifications are required.

First, given a sequence of homomorphisms $\{ h:FA_{\ell} \to FS(U)\}$, that are associated with $\ell$-tuples, we can apply the techniques that were used
in the proof of theorem 1.1, and extract a subsequence that converges and factors through a  resolution: $L_1 \to L_2 \to  \ldots \to L_f$, where
$L_f$ is a limit group and the sequence of shortened homomorphisms that converges to $L_f$, is obtained from a sequence of
homomorphisms: $\{ \tilde h :FS_{\ell} \to FS(U)\}$ that have uniformly bounded images. i.e., there exists some constant $c$, so that the images of a fixed
generated set of $FS_{\ell}$ is mapped by each of the homomorphisms to $\ell$ elements in $FS(U)$ of length bounded by $c$.

Note that since $FS(U)$ is countably generated, there are countably many elements of length bounded by $c$ in $FS(U)$. At this point we change the metric on
$FS(U)$, and give each generator $u \in U$, a length which is its degree in $FA_k$. If the union of the images of the homomorphisms $\{\tilde h\}$ is not a finite
set of $\ell$-tuples in $FS(U)$, then from the sequence $\{\tilde h\}$ it is possible to extract a subsequence that converges into an action of
the terminal limit group $L_f$ on a simplicial tree with trivial edge groups. This means that the constructed resolution is not single ended.

Hence, if the constructed resolution is single ended, then the shortened homomorphisms into the terminal limit group $L_f$, $\{\tilde h\}$, have only finitely
many possible images, so after passing to a further subsequence we can assume that the terminal homomorphisms, $\{\tilde h\}$, are in fact constant
homomorphisms.

Therefore, from any sequence of homomorphisms: $\{h:FS_{ell} \to FS(U)\}$ it is possible to extract a subsequence that converges into a finite
resolution (as in theorem 1.1), and if the resolution is single ended, the resolution can be encoded by a finite sequence. Given a set of $\ell$-tuples in $FA_k$,
we consider the countably many resolutions that are obtained from  convergent subsequences, that are extracted from the homomorphisms,
$h:FS_{\ell} \to FS(U)$, that are associated with the given set of $\ell$-tuples. We assume that all these resolutions are single ended.

The countable set of resolutions can be naturally ordered, and by the compactness argument that is used to prove theorem 1.2, the sequence
of resolutions that are extracted from convergent sequences,  can be replaced by a finite MR diagram that satisfy the conclusions of the theorem 1.2.

Once there exists a single ended MR diagram that is associated with a set of $\ell$-tuples from $FA_k$, we can replace the MR diagram with a minimal
complexity one, and claim that the minimal complexity MR diagram is essentially unique up to isomorphism following the arguments that were used in the
proof of theorem 1.16.

\smallskip
So far we have assumed that the algebra $FA_k$ is defined over a finite field $F$. This enabled us to deduce that a given subsequence of
homomorphisms that converge into a single ended resolution terminates with a
fixed homomorphism. If $F$ is a countable field, we can still deduce that a single ended resolution terminates with homomorphisms with bounded image. i.e.,
that the shortened homomorphisms that converge to terminal the limit group $L_f$ in a resolution map a fixed set of generators of $FS_{\ell}$ onto
homogeneous elements
of uniformly bounded degree in $FA_k$. This implies that these terminal shortened homomorphisms satisfy a system of equations over $FA_k$ that translates
to a system of equations over the countable field $F$.

Therefore, over a countable field $F$, a resolution can also be encoded with a finite sequence, and the compactness argument that was used in theorem 1.2 and for
algebras $FA_k$ over a finite field can still be applied. Having a single ended MR diagram we apply the arguments that prove theorem 1.16, and obtain a unique
minimal complexity single ended MR diagram for a set of tuples over a free associative algebra that is defined over a countable field.

\smallskip
Given a system of equations over a free algebra,
we can look at its associated variety, i.e., its set of solutions. These are in general not tuples
of homogeneous elements but rather of non-homogeneous ones. With the tuples of non-homogeneous elements that form the variety
(the solution to a system of equations) we associate their top homogeneous tuples. With these top homogeneous tuples we can associate an MR diagram
and if its single ended, a canonical (minimal complexity) diagram.

In a forthcoming work with A. Atkarskaya [At-Se] we prove that in the case of a single ended diagram, this canonical diagram, that is built from the top homogeneous
parts of the set of solutions, gives a description of the set of solutions, i.e., of the set of non-homogeneous solutions to a given system of equations.

\smallskip
Let $LA_k$ be a free Lie algebra that is generated by $k$ elements over a finite field. $LA_k$ can be embedded into $FA_k$ by the standard embedding that maps
$[x,y]$ to $xy-yx$. Hence, every with an  ordered $\ell$-tuple in $LA_k$  we can associate an ordered $\ell$-tuple in $FA_k$. Therefore, with a variety over
$LA_k$ we can associate a set of $\ell$-tuples in $FA_k$.

With the set of $\ell$-tuples in $FA_k$ we associate their top homogeneous parts. With the set of their top homogeneous parts we can associate an MR diagram,
and if the MR diagram is single ended, we associate the canonical minimal complexity MR diagram that is associated with the set of homogeneous tuples.
Now, the work on varieties in a free associative algebra gives a description of the non-homogeneous elements in the image of the variety over $LA_k$ in $FA_k$
using the structure of the canonical MR diagram that is associated with the top homogeneous elements of the non-homogeneous elements in the image of the variety.


\vskip 2em


\end
\smallskip
\Refs

\widestnumber\key{XX-XXX}

\ref\key At-Se
\by A. Atkarskaya and Z. Sela
\paper Monomial varieties over a free associative algebra
\paperinfo in preparation
\endref
\vskip 2em

\ref\key Be
\by G. Berk
\paper Canonicality of Makanin-Razborov diagrams - counterexample
\jour Annals Inst. Fourier (Grenoble)
\vol 70 \yr 2020
\pages 2027-2047
\endref
\vskip 2em



\ref\key Gu
\by V. Guirardel
\paper Actions of finitely generated groups on R-trees
\jour Annals Inst. Fourier (Grenoble)
\vol 58 \yr 2008
\pages 159-211
\endref
\vskip 2em


\ref\key Ja-Se
\by E. Jaligot and Z. Sela
\paper Makanin-Razborov diagrams over free products
\jour Illinois jour. of math.
\vol 54 \yr 2010 \pages 19-68
\endref
\vskip 2em

\ref\key Je
\by A. Jez
\paper Recompression: A simple and powerful technique for word equations
\jour Jour. of the ACM
\vol 63 \yr 2016 \pages 151-
\endref
\vskip 2em


\ref\key  Ma
\by G. S. Makanin
\paper Equations in a free group
\jour Math. USSR Izvestiya
\vol 21 \yr 1983
\pages 449-469
\endref
\vskip 2em





\ref\key Po1
\by B. Poizat
\paper Groupes stables avec types generiques reguliers
\jour Journal of symbolic logic
\vol 48  \yr 1983 \pages 339-355
\endref
\vskip 2em

\ref\key Po2
\bysame
\paper A course in model theory
\paperinfo Springer, 2000
\endref
\vskip 2em

\ref\key Ra1
\by A. A. Razborov
\paper On systems of equations in a free group
\jour Math. USSR Izvestiya
\vol 25 \yr 1985 \pages 115-162
\endref
\vskip 2em

\ref\key Ra2
\bysame
\paper On systems of equations in a free group
\paperinfo Ph.D. thesis, Steklov Math. institute, 1987
\endref
\vskip 2em


\ref\key Re-We
\by C. Reinfeldt and R. Weidmann
\paper Makanin-Razborov diagrams for hyperbolic groups
\jour Annales Math. Blaise Pascal
\vol 26 \yr 2019 \pages 119-208
\endref
\vskip 2em



\ref\key Se1
\by Z. Sela
\paper Diophantine geometry over groups I: Makanin-Razborov diagrams
\jour Publications Mathematique de l'IHES
\vol 93 \yr 2001 \pages 31-105
\endref
\vskip 2em

\ref\key Se2
\bysame
\paper Diophantine geometry over groups II: Completions, closures and formal solutions
\jour Israel jour. of Math.
\vol 134 \yr 2003 \pages 173-254
\endref
\vskip 2em


\ref\key Se3
\bysame
\paper Diophantine geometry over groups IV: An iterative procedure for validation of a sentence
\jour Israel jour. of Math.
\vol 143 \yr 2004 \pages 1-130
\endref
\vskip 2em

\ref\key Se4
\bysame
\paper Diophantine geometry over groups V$_1$: Quantifier elimination I
\jour Israel jour. of Mathematics
\vol 150 \yr 2005 \pages 1-197
\endref
\vskip 2em


\ref\key Se5
\bysame
\paper Diophantine geometry over groups V$_2$: Quantifier elimination II
\jour GAFA
\vol 16 \yr 2006 \pages 537-706
\endref
\vskip 2em


\ref\key Se6
\bysame
\paper Diophantine geometry over groups VII: The elementary theory of a
hyperbolic group
\jour Proceedings of the LMS
\vol 99 \yr 2009 \pages 217-273
\endref
\vskip 2em

\ref\key Se7
\bysame
\paper Diophantine geometry over groups VIII: Stability
\jour Annals of Math.
\vol 177 \yr 2013 \pages 787-868
\endref
\vskip 2em


\ref\key Se8
\bysame
\paper Diophantine geometry over groups IX: Envelopes and imaginaries
\paperinfo preprint
\endref
\vskip 2em

\ref\key Se9
\bysame
\paper Word equations I: Pairs and their Makanin-Razborov diagram
\paperinfo preprint
\endref
\vskip 2em

\ref\key Te-Zi
\by K. Tent and M. Ziegler
\paper A course in model theory
\paperinfo Lecture notes in logic, ASL, Cambridge  university press 2012
\endref
\vskip 2em



\ref\key We
\by R. Weidmann
\paper On accessibility of finitely generated groups
\jour Quarterly journal of math.
\vol 63 \yr 2012 \pages 211-225
\endref
